\documentclass{tac}

\usepackage{modalops}
\usepackage{hyperref}
\usepackage{amsfonts}
\usepackage{tikz-cd}
\usepackage{amsmath,amssymb}
\usepackage{dirtytalk}
\usepackage{comment}
\usepackage[capitalise,nameinlink]{cleveref}
\usepackage{longtable}

\newcommand{\up}[1]{#1^{\uparrow}}
\newcommand{\fork}[1]{(#1 [1])^+}

\newcommand{\N}{\mathbb{N}}
\newcommand{\Z}{\mathbb{Z}}

\newcommand{\subred}{\begin{tikzcd}[ampersand replacement=\&]
	{} \& {}
	\arrow[from=1-1, to=1-2, {Circle[open]}->, start anchor={[xshift=-1.245ex]}, end anchor={[xshift=-3.245ex]}]
\end{tikzcd} \hspace{-1.955em}}

\newcommand{\Set}{\textbf{\textup{Set}}}
\newcommand{\KFr}{\textbf{\textup{KFr}}}
\newcommand{\MA}{\textbf{\textup{MA}}}
\newcommand{\HA}{\textbf{\textup{HA}}}
\newcommand{\CAMA}{\textbf{\textup{CAMA}}}
\newcommand{\CABA}{\textbf{\textup{CABA}}}
\newcommand{\DGrph}{\textbf{\textup{DGrph}}}

\newcommand{\lf}{\textit{lf}}
\newcommand{\fin}{\textit{fin}}

\newcommand{\op}{\textup{op}}
\newcommand{\id}{\textup{id}}

\DeclareMathOperator{\Pro}{Pro-}

\DeclareMathOperator{\Lex}{Lex}
\DeclareMathOperator{\Refl}{Pos}
\DeclareMathOperator{\depth}{d}
\DeclareMathOperator{\width}{wt}
\newcommand{\exter}{\delta^e}
\newcommand{\inter}{\delta^i}
\DeclareMathOperator{\numberext}{e}
\DeclareMathOperator{\Image}{Im}
\DeclareMathOperator{\NExt}{NExt}
\DeclareMathOperator{\Ext}{Ext}
\DeclareMathOperator{\Chains}{BC}

\title[Profinite interior algebras]{Finite coproducts, coregularity and coexactness for profinite interior algebras}

\author{Matteo De Berardinis}
\address{University of Salerno, University of Amsterdam}
\eaddress{mdeberardinis@unisa.it}

\keywords{profiniteness, interior algebras, regularity, Barr exactness}

\amsclass{03B45, 03G05, 18A30, 18A32, 18C35, 18E08}

\begin{document}
\maketitle
\begin{abstract}
In previous articles, we showed that the category of profinite \textit{L}-algebras (where \textit{L} is a normal modal logic with the finite model property) is monadic over \textbf{Set}. Then, we developed sequent calculi for extensions of the language of \textit{L} with infinitary conjunctions and disjunctions, proving completeness with respect to profinite \textit{L}-algebras and relating syntactic properties of the calculi with regularity/exactness properties of the category opposite to profinite \textit{L}-algebras. In this paper, we focus on the algebraic perspective: we characterize those \textit{L} extending \textsf{S4} whose profinite algebras enjoy such categorical properties.
\end{abstract}

\section{Introduction}
In \cite{GdB24} we proved that the category $\Pro L\MA_\fin$ of profinite $L$-algebras, for any normal modal logic $L$ with the finite model property, is monadic over $\Set$. This result ensures that it is possible to describe profinite $L$-algebras in terms of (possibly infinitary) operations and equations, while guaranteeing the existence of the free objects. We then developed infinitary sequent calculi (see \cite{infinitary}) for propositional modal languages with infinitary conjunctions and disjunctions, enriching Maehara-Takeuti's infinitary extension of the sequent calculus \textbf{LK} (introduced in \cite[Chapter 4]{Tak}) with peculiar modal rules, while keeping well-foundedness of the proofs. A Lindenbaum algebra construction proves that profinite $L$-algebras provide a complete semantics for such calculi. 

The two aforementioned descriptions of profinite $L$-algebras are then compared in \cite{infinitary}: bridge theorems relating categorical properties of $\Pro L\MA_\fin$, or of its dual category $L\KFr_\lf$ of \emph{locally finite} Kripke frames for $L$, and syntactic properties of the corresponding calculi are established. Namely, as every monadic category over $\Set$, $\Pro L\MA_\fin$ (hence its dual $L\KFr_\lf$) is both complete and cocomplete and it has both a regular factorization system and a coregular one.

On the one hand, it turns out that the coincidence of the two factorization systems (which implies that both all epis and all monos are regular --- actually, it is equivalent to the fact that all epis are regular in $L\KFr_\lf$) corresponds to a strong version of Beth's Definability Property. On the other hand, the fact that all monos are regular in $L\KFr_\lf$ is equivalent to an infinitary version of the standard Beth's Property. Stability of the coregular factorization system in $L\KFr_\lf$ (in \cite{GZ} called --- together with the existence of finite limits --- \emph{r-regularity}, as opposed to the standard regularity, which takes into account the regular factorization system) is equivalent to amalgamability for profinite $L$-algebras and it corresponds to Craig's Interpolation Theorem for the global consequence relation. If, in addition to r-regularity, we require that equivalence relations are in one-to-one correspondence with quotients, we obtain \emph{r-Barr exactness},
%rev . r-Barr exactness is also considered in \cite{infinitary}.
a notion also considered in \cite{infinitary}. In many cases, the aforementioned results are similar to those existing for the analogous analysis in the customary finitary fragment of our language (see \cite{GZ}), with the difference that profinite algebras have to be replaced with finitely presented ones.

This paper can be seen as a continuation of \cite{infinitary}. When all monomorphisms in $L\KFr_\lf$ are regular --- for instance, whenever $L$ extends \textsf{K4} (meaning that the Kripke models of $L$ are transitive) --- then the notions of r-regularity and r-Barr exactness coincide with their standard (non r-) versions. We focus on extensions of the logic \textsf{S4}: the Kripke models are preorders (the relation is reflexive, as well as transitive); the algebraic models are called \emph{interior} algebras, because their axioms encode the properties defining the topological interior. In particular, we characterize the logics $L$ extending \textsf{S4} for which $L\KFr_\lf \simeq (\Pro L\MA_\fin)^\op$ is a regular/Barr exact category (we leave the full characterization of the extensions of \textsf{K4} to future work). For this classification, we present the logics semantically rather than axiomatically: the assumption of the finite model property ensures that a logic is completely determined by the class of its finite models.

The paper is structured as follows. In Section~\ref{sec:profinite}, we recall the basic definitions regarding algebraic and Kripke semantics for normal modal logics (see \cite{misha} for a complete overview on modal logic), as well as some results from \cite{GdB24}.

In Section~\ref{sec:colimits}, we give an explicit description of the finite limits in $L\KFr_\lf$. In particular, the construction of the finite products will turn out to be strongly connected to the so-called \emph{universal model} construction, namely a construction widely investigated in the modal logic literature for transitive modal systems \cite{leo,valentin,be1,be2,ungh}, in connection to problems like representation theorems, atomicity of finitely generated free algebras, local finiteness of subvarieties, etc.

In Section~\ref{sec:AP}, we fully characterize those normal modal logics $L \supseteq \textsf{S4}$ with the finite model property for which $L\KFr_\lf$ is a regular category, i.e. for which $\Pro L\MA_\fin$ has the amalgamation property (see \cite{infinitary}). The resulting classification follows Maksimova’s parallel results for the finitary fragment \cite{Maksimova1979InterpolationTI,Maksimova1980InterpolationTI} (with some notable exceptions taken from \cite{Maksimova1982AbsenceOT}).

Finally, we characterize Barr exactness. On the one hand, the attempt of Section~\ref{sec:coexactness} to describe equivalence relations leads us to the formulation of necessary conditions for $L\KFr_\lf$ to be Barr exact, by means of the absence of certain configurations among the Kripke models of $L$. On the other hand, in Section~\ref{sec:presheaves} we present an adjunction $\widehat{\mathsf{C}} \longrightarrow \KFr$ for the category of presheaves over any small category $\mathsf{C}$, which restricts to an equivalence with $L\KFr_\lf$, for some trivial $L$'s. Since the category of presheaves over a small category is Barr exact (being a topos) it turns out that there are exactly five logics $L$ above \textsf{S4} with the finite model property such that $L\KFr_\lf$ is a Barr exact category.

\section{Profinite modal algebras}\label{sec:profinite}
In this Section we review (and integrate) some notions and results from \cite{GdB24} in order to supply the  framework of the paper.

Throughout the paper, $L$ will denote a normal modal logic with the finite model property. A \emph{normal modal logic} is a set $L$ of modal formulas containing all tautologies of classical propositional calculus and the \emph{distribution axiom} $\nec (x \to y) \to (\nec x \to \nec y)$ and closed under \emph{uniform substitution} and the rules $\varphi \to \psi, \varphi/\psi$ (\emph{modus ponens}) and $\varphi/\nec \varphi$ (\emph{necessitation}). The minimum normal modal logic is denoted by \textsf{K}. The family of normal modal logics containing (we sometimes say extending or above) a logic $L$ is denoted by $\NExt L$. If $\Gamma$ is a set of modal formulas, $L + \Gamma$ denotes the smallest $L' \in \NExt L$ containing $\Gamma$. A normal modal logic is said to have the \emph{finite model property} if, for each $\varphi \notin L$, there exists a finite (algebraic) model of $L$ that falsifies $\varphi$.

The algebraic semantics for the minimum normal modal logic \textsf{K} \cite{misha} is given by \emph{modal algebras}: a modal algebra $(B,\pos)$ is a Boolean algebra $B$ endowed with a finite-join preserving operator $\pos \colon B \longrightarrow B$ such that
$$\pos (x \vee y) = \pos x \vee \pos y ~~~~~\text{and}~~~~~ \pos \bot = \bot.$$
The operator $\pos$ is called the ‘possibility’ operator and its dual ($\nec x \coloneqq \neg \pos \neg x$) is called the ‘necessity’ operator. Modal algebras and Boolean morphisms preserving $\pos$ form the category $\MA$. We will consider varieties of modal algebras generated by their finite members. Such varieties (called `finitely approximable') are in bijective correspondences with normal modal logics with the finite model property \cite{misha}: if $L$ is such a logic, we call $L$\emph{-algebras} their algebraic models and $L\MA$ the corresponding variety, viewed as a category.

\emph{Profinite $L$-algebras} are defined as the \emph{Pro-Completion} (= the categorical construction of freely adjoining cofiltered limits to a category, see \cite{Stone}) of the full subcategory $L\MA_\fin$ of $L\MA$ given by finite $L$-algebras.
Below, we shall recover them as a full subcategory of $\CAMA_\infty$, the category of \emph{complete}, \emph{atomic} and \emph{completely additive} (= $\pos$ commutes with arbitrary joins) modal algebras and Boolean morphisms preserving arbitrary joins (hence arbitrary meets), as well as $\pos$. To this aim, we first recall Goldblatt-Thomason duality of $\CAMA_\infty$ with Kripke frames.

A \emph{Kripke frame} $(W,\prec)$ is a set $W$ endowed with a binary relation $\prec$ on it (most of the times we will omit to indicate the binary relation). A \emph{p-morphism} from a Kripke frame $W$ to a Kripke frame $V$ is a function $f \colon W \longrightarrow V$ with the following properties:

\begin{enumerate}
    \item (stability) for $w, w' \in W$, $w \prec w' \implies f(w) \prec f(w')$;
    \item (openness) for $v' \in V$, $w \in W$, $f(w) \prec v' \implies \exists w' \in W \text{ such that } w \prec w' ~\&~ f(w') = v'$.
\end{enumerate}
Kripke frames and p-morphisms form the category $\KFr$. We call $\DGrph$ the category of Kripke frames and stable functions between them. Extending Tarski duality between $\CABA$ (complete atomic Boolean algebras) and $\Set$, we obtain the following.
\theorem[Goldblatt-Thomason duality \cite{thomason}]\label{teo:duality}
$\CAMA_\infty$ is dual to $\KFr$.
\endtheorem
\proof(Sketch)
On the modal algebras side, the duality functors make the set of atoms of a complete, atomic and completely additive modal algebra into a Kripke frame, by setting $a \prec b$ iff $a \leq \pos b$. On the other side, for a Kripke frame $W$, the Boolean algebra $\mathcal{P}(W)$ is turned into a modal algebra by setting $\pos X := \{w \in W\ \vert\ \exists w' (w \prec w' ~\&~ w' \in X)\}$.
\endproof

Given a logic $L$, we can consider the full subcategory $L\CAMA_\infty$ formed by the complete, atomic and completely additive modal algebras belonging to $L\MA$; its dual is the category $L\KFr$, whose objects are precisely those Kripke frames validating $L$. Since finite modal algebras are complete, atomic and completely additive, the category of finite modal $L$-algebras $L\MA_\fin$ is dual to $L\KFr_\fin$, the full subcategory of $L\KFr$ consisting of its finite Kripke frames. In order to study the Pro-Completion of $L\MA_\fin$, we can describe the dual construction, called \emph{Ind-Completion}, over the category $L\KFr_\fin$.

\definition\label{defn:lf}
A Kripke frame $W$ is said to be \emph{locally-finite} if, for all $w \in W$,
$$w^* \coloneqq \{w' \in W\ \vert\ \exists w_0,\dots,w_n \in W \text{ such that } w=w_0 \prec \dots \prec w_n=w'\}$$
is a finite set.
\enddefinition
We denote by $L\KFr_\lf$ the full subcategory of $L\KFr$ consisting of its locally-finite Kripke frames.

\proposition\label{prop:colim}
$\KFr$ is cocomplete; colimits are inherited by the full subcategory $L\KFr_\lf$ and created by the forgetful functor to $\Set$.
\endproposition
\proof
See \cite{GdB24}.
\endproof

\theorem\label{teo:ind}
$L\KFr_\lf$ is the Ind-Completion of $L\KFr_\fin$ via the full embedding.%$L\KFr_\fin \longrightarrow L\KFr_\lf$.
\endtheorem
\proof
See \cite{GdB24}.
\endproof

\corollary\label{cor:lim}
$L\KFr_\lf$ has all limits.
\endcorollary
\proof
See \cite{GdB24}.
\endproof

\theorem\label{teo:produality} The Goldblatt-Thomason duality restricts to a duality between the category of profinite modal $L$-algebras $\Pro L\MA_\fin$ (the Pro-Completion of $L\MA_\fin$) and the category  $L\KFr_\lf$.
\endtheorem

From the above picture, it follows that $\Pro L\MA_\fin$ is complete and cocomplete; the forgetful functor from $L\CAMA_\infty$ to $\Set$ preserves limits, but does not have an adjoint; on the contrary, its restriction to $\Pro L\MA_\fin$ does have an adjoint.

\theorem\label{teo:mon}
The forgetful functor $\Pro L\MA_\fin \longrightarrow \Set$ is monadic.
\endtheorem
\proof
This is proved in \cite{GdB24} in two ways, either by a direct application of Beck theorem \cite{CWM}, or by combining Manes monadicity theorem~\cite{manes} with specific facts about quotients of profinite algebras~\cite{gehrke}.
\endproof

Notable examples of normal modal logics with the finite model property are the logics \textsf{K4}, \textsf{S4}, \textsf{GL}, \textsf{Grz}, \textsf{GL.3}, \textsf{Grz.3}, etc: the Table below\footnote{In the Table $\nec^+ x$ stands for $x \wedge \nec x$; a relation $\prec$ is said to be \emph{locally linear} iff it satisfies the condition $(w\prec v_1 ~\&~ w\prec v_2) \Rightarrow ( v_1\prec v_2 \vee v_1=v_2 \vee v_2\prec v_1)$ for all $v,w_1,w_2$.}
reports the axiomatization of such logics and the characterization of their (locally) finite Kripke frames (see \cite{misha} for more information and for the relevant proofs).

\begin{center}
\begin{longtable}{||l|l|l||} \hline

Logic  $L$ & Axioms & Corresponding (locally) \\
& & finite Kripke frames \\ \hline

%rev $K$ & ~~~~~~~ - &~~~~~~ - \\ \hline

\textsf{K4} & \textsf{K} + $\nec x \to \nec \nec x$ & $\prec$ is transitive \\ \hline

\textsf{S4} & \textsf{K4} + $\nec x \to x$ & $\prec$ is reflexive and \\
& & transitive \\ \hline

\textsf{S4.2} & \textsf{S4} + $\pos \nec x \to \nec \pos x$ & $\prec$ is reflexive, transitive \\
& & and locally confluent \\ \hline

\textsf{S4.3} & \textsf{S4} + & $\prec$ is reflexive, transitive \\
& $\nec (\nec x_1 \to x_2) \vee \nec (\nec x_2 \to x_1)$ & and locally linear \\ \hline

\textsf{Grz} & \textsf{S4} + & $\prec$ is a partial order \\
& $\nec (\nec (x \to \nec x) \to x) \to x$ & \\ \hline

\textsf{Grz.3} & \textsf{S4.3} + \textsf{Grz} & $\prec$ is a locally linear \\
& & partial order \\ \hline

\textsf{S5} & \textsf{S4} + $x \to \nec \pos x$ & $\prec$ is an equivalence \\
& & relation \\ \hline

\textsf{GL} & \textsf{K4} + $\nec (\nec x \to x) \to \nec x$ & $\prec$ is irreflexive and \\
& & transitive \\ \hline

$\textsf{GL.3}$ & \textsf{GL} + & $\prec$ is irreflexive, \\
& $\nec^+ (\nec^+ x_1 \to x_2) \vee \nec^+ (\nec^+ x_2 \to x_1)$ & transitive and \\
& & locally linear \\
\hline
\end{longtable}
\end{center}

\section{Finite coproducts}\label{sec:colimits}
Being monadic, the forgetful functor $\Pro L\MA_\fin \longrightarrow \Set$ preserves limits; as a consequence, limits in $\Pro L\MA_\fin$ are computed as in $\Set$.

Colimits in $\Pro L\MA_\fin$ (dually, limits in $L\KFr_\lf$, by Theorem~\ref{teo:produality}) are more complicated. The idea to compute products in $L\KFr_\lf$ as in the underlying category $\DGrph$ (Kripke frames and stable functions) is wrong: under appropriate assumptions, the projections from the $\DGrph$-product\footnote{The $\DGrph$-product of two Kripke frames $W_0$ and $W_1$ has the cartesian product of $W_0$ and $W_1$ as underlying set and it is equipped with the product relation.}
are p-morphisms; however, the morphism induced by the universal property is not necessarily open. A consequence of Theorem~\ref{teo:ind} (together with the theory of Ind-Completions --- see \cite{Stone}) is that $L\KFr_\lf$ is equivalent to the category $\Lex(L\MA_\fin,\Set)$ of finite-limits-preserving functors $L\MA_\fin \longrightarrow \Set$, and natural transformations between them. Limits in $\Lex(L\MA_\fin,\Set)$ can be computed point-wise; however, the equivalence $L\KFr_\lf \simeq \Lex(L\MA_\fin,\Set)$ is rather involved.

In this section, we explicitly describe finite limits in $L\KFr_\lf$, under the assumption that $L \in \NExt \textsf{S4}$ (so that the objects of $L\KFr_\lf$ are locally finite transitive and reflexive Kripke frames, i.e.\ locally finite preorders); finite colimits in $\Pro L\MA_\fin$ are then obtained via the restricted Goldblatt-Thomason duality.

We introduce some notions that will be useful in this and in the following sections. For a Kripke frame $W$ and $w \in W$, we set $\up{w} \coloneqq \{w' \in W\ \vert\ w \prec w'\}$.
\definition
Given a Kripke frame $W$, a subset $G \subseteq W$ is said to be a \emph{generated subframe} of $W$ if, for all $w \in W$,
$$w \in G ~~~~~\implies~~~~~ \up{w} \subseteq G.$$
\enddefinition
A generated subframe $G$ of $W$ can be seen as a Kripke frame, by restricting the binary relation of $W$ to $G$,
%rev $G$ generated subframe of $W$ can be seen as a Kripke frame, by taking the restriction of the binary relation on $W$,
so that the inclusion $G \longrightarrow W$ becomes a p-morphism. Observe that, if $W$ is locally finite, then $G$ is locally finite, too. Moreover, the logic of $G$ contains the logic of $W$ (see \cite{misha}). As a consequence, if $W \in L\KFr_\lf$, then $G \in L\KFr_\lf$. Recall that, in Definition~\ref{defn:lf}, we introduced the set
$$w^* \coloneqq \{w' \in W\ \vert\ \exists w_0,\dots,w_n \in W \text{ such that } w=w_0 \prec \dots \prec w_n=w'\}$$
for $w \in W$. Observe that $w^*$ is the smallest generated subframe of $W$ containing $w$. Moreover, $w^*$ coincides with $\up{w}$ if and only if $W$ is reflexive and transitive.

For a locally finite transitive Kripke frame $W$ (reflexivity is not required here) a notion of depth can be defined, allowing for powerful induction arguments in proofs.
%rev it is defined a notion of depth, leading to powerful induction arguments in proofs.
The \emph{depth} $\depth(w)$ of a point $w \in W$ is defined as the maximum cardinality of sequences
\begin{align*}
    w = w_1 \prec w_2 \prec \dots \prec w_n
\end{align*}
in $W$, such that $w_i \neq w_{i+1}$ and $w_{i+1} \nprec w_i$ for $i = 1, \dots, n-1$ (such sequences are called \emph{chains}). Since $W$ is locally finite and transitive, $\depth(w)$ is a positive natural number.
%rev $\depth(w)$ is a natural number greater or equal than $1$, $W$ being locally finite and transitive.
We can then define $\depth(W)$ as the maximum of the depths of the points in $W$, with the convention that $\depth(W) = 0$ if $W$ is empty and $\depth(W) = \omega$ if the maximum does not exist.

A transitive Kripke frame can be partitioned into its irreflexive points and its clusters. A \emph{cluster} in a transitive Kripke frame $W$ is a subset of $W$ maximal (with respect to inclusion) among the nonempty subsets $C \subseteq W$ such that the restriction of the binary relation to $C$ is maximal, i.e., for each $w, w' \in W$,
\begin{align*}
w, w' \in C ~~~~~\implies~~~~~ w \prec w' ~\&~ w' \prec w.
\end{align*}
Observe that the depth is constant on clusters, so that, given a cluster $C \subseteq W$, the non-zero natural number $\depth(C)$ is well defined.
%rev it is well defined the non-zero natural number $\depth(C)$.
The cluster $C$ is said to be \emph{external} if $\depth(C) = 1$, \emph{internal} otherwise.

Fix $L \in \NExt \textsf{S4}$ with the finite model property; it is straightforward to prove that the terminal object in $L\KFr_\lf$ is the one-element reflexive Kripke frame.

We now describe equalizers in $L\KFr_\lf$. Consider a pair of parallel p-morphisms $f, g \colon W \longrightarrow V$ in $L\KFr_\lf$; set
\begin{align*}
E_{fg} \coloneqq \{w \in W\ \vert\ \forall w' (w' \in w^* \implies f(w') = g(w'))\}.
\end{align*}
Observe that $E_{fg}$ is the biggest generated subframe of $W$ contained in the $\Set$-equalizer $\{w \in W\ \vert\ f(w)=g(w)\} \subseteq W$.
\lemma\label{lem:eq}
The inclusion $E_{fg} \longrightarrow W$ is the equalizer of $f$ and $g$ in $L\KFr_\lf$.\footnote{The assumption $L \in \NExt \textsf{S4}$ is not required here.}
\endlemma
\proof
Consider another p-morphism $h \colon A \longrightarrow W$ such that $f h = g h$. $\Image(h)$ is a generated subframe of $W$ and it is contained in $\{w \in W\ \vert\ f(w)=g(w)\}$; as a consequence, $\Image(h) \subseteq E_{fg}$.
\endproof

Finally, we describe binary products in $L\KFr_\lf$. Let us start with $L = \textsf{S4}$. Consider two locally finite preorders $W_0$ and $W_1$. For each $n \geq 0$, define a cone
\[\begin{tikzcd}\label{eq:prodn}
	{X^n} & {W_0} \\
	& {W_1}
	\arrow["{p_0^n}", from=1-1, to=1-2]
	\arrow["{p_1^n}"', from=1-1, to=2-2]
\end{tikzcd}\tag{$P_n$}\]
in the following way. $X^0$ is the initial object $\emptyset$ and $p_0^0$ and $p_1^0$ are the unique morphisms induced by the universal property (the empty functions). Suppose that $X^n$, $p_0^n$ and $p_1^n$ have been defined. Put $X^{n,+}$ equal to the set of pairs $(Y,G)$ such that (i) $Y$ is a non-empty finite subset of the $\Set$-product $W_0 \times W_1$ and $G$ is a finite generated subframe of $X^n$ comprising at least a point of depth $n$, if $n > 0$; (ii) $\forall y \in Y$, $\up{\pi_i(y)} = \pi_i(Y) \cup p_i^n(G)$, for $i = 0, 1$;\footnote{$\pi_i \colon W_0 \times W_1 \longrightarrow W_i$ are the projections of the cartesian product in \Set.} (iii) $Y \nsubseteq \{(p_0^n(x),p_1^n(x))\ \vert\ x \in X^n \text{ s.t. } \up{x} = G\}$. Set
\begin{align*}
X^{n+1} \coloneqq X^n \cup \{(y,Y,G)\ \vert\ (Y,G) \in X^{n,+} ~\&~ y \in Y\}
\end{align*}
(assuming the union to be disjoint) and, for $x, x' \in X^{n+1}$,
\begin{align*}
x \prec x' \iff \begin{cases}
    x, x' \in X^n \text{ with } x \prec x', \qquad \text{or} \\
    x = (y,Y,G) \text{ for } (Y,G) \in X^{n,+} ~\&~ y \in Y,\, x' \in G, \qquad \text{or} \\
    x = (y,Y,G),\, x' = (y',Y,G) \text{ for } (Y,G) \in X^{n,+} ~\&~ y, y' \in Y.
\end{cases}
\end{align*}
Moreover, set
\begin{align*}
p_i^{n+1}(x) \coloneqq  \begin{cases}
    p_i^n(x)       & \quad \text{if } x \in X^n \\
    \pi_i(y)       & \quad \text{if } x=(y,Y,G) \text{ for } (Y,G) \in X^{n,+} ~\&~ y \in Y
\end{cases}
\end{align*}
for $i = 0, 1$. Observe that $X^n \subseteq X^{n+1}$ is a generated subframe and that $p_i^{n+1}|_{X^n} = p_i^n$ for $i = 0, 1$. If we set $X \coloneqq \cup_n X^n$ and $p_i \coloneqq \cup_n p_i^n \colon X \longrightarrow W_i$ for $i = 0, 1$, we obtain the following diagram
\[\begin{tikzcd}\label{eq:prod}
	{X} & {W_0} \\
	& {W_1}
	\arrow["{p_0}", from=1-1, to=1-2]
	\arrow["{p_1}"', from=1-1, to=2-2]
\end{tikzcd}\tag{$P$}\]

\remark
We compute the first stages of the construction above, instantiated on
\begin{center}
\begin{tikzpicture}[
  vertex/.style = {circle, draw,
  minimum size=26, inner sep=0,
  fill=white},
  cluster/.style = {circle, draw,
  minimum size=22mm, inner sep=0,
  fill=white},
  vertex1/.style = {vertex, fill=red!20!white},
  vertex2/.style = {vertex, fill=orange!30!white},
  vertex3/.style = {vertex, fill=blue!30!white},
  vertex4/.style = {vertex, fill=teal!30!white},
]

  \draw[thick]
    (-5,0) node {$W_0=$}
    (-3,0.7) node[cluster] (A) {}
    ([shift=({180:0.5})]A) node[vertex] {$1$}
    ([shift=({0:0.5})]A) node[vertex] {$2$}
    (-3,-1.3) node[vertex] (B) {$\rho$}
    (A) -- (B);

  \draw[thick]
    (1,0) node {$W_1=$}
    (2,1) node[vertex] (1') {$a$}
    (4,1) node[vertex] (2') {$b$}
    (3,-1) node[vertex] (r') {$\delta$}
    (r') -- (1')
    (r') -- (2');
\end{tikzpicture}
\end{center}
A pair $(Y,G) \in X^{0,+}$ is such that $G=\emptyset$ and $Y \subseteq W_0 \times W_1$, with $\pi_i(y)^\uparrow = \pi_i(Y)$ for each $y \in Y$ and $i = 0,1$ (condition (iii) is obviously verified, since $X^0 = \emptyset$). We cannot have $y \in Y$ with $\pi_0(y) = \rho$: otherwise, $\{\rho,1,2\} = \rho^\uparrow = \pi_0(y)^\uparrow = \pi_0(Y)$, hence there exists $y' \in Y$ such that $\pi_0(y') = 1$; but then $\{\rho,1,2\} = \pi_0(Y) = \pi_0(y')^\uparrow = 1^\uparrow = \{1,2\}$. Similarly, we cannot have $y \in Y$ such that $\pi_1(y) = \delta$. With a similar argument, it is possible to exclude the sets of pairs $\{(1,a)\}$, $\{(1,b)\}$, $\{(1,a),(1,b)\}$, and so on. The frame $X^1$, containing the clusters of depth $1$ of the product, will then be
\begin{center}
\begin{tikzpicture}[
  vertex/.style = {circle, draw,
  minimum size=26, inner sep=0,
  fill=white},
  cluster/.style = {circle, draw,
  minimum size=22mm, inner sep=0,
  fill=white},
  vertex1/.style = {vertex, fill=red!20!white},
  vertex2/.style = {vertex, fill=orange!30!white},
  vertex3/.style = {vertex, fill=blue!30!white},
  vertex4/.style = {vertex, fill=teal!30!white},
]

  \draw[thick]
    (-2,0) node[cluster] (A) {}
    ([shift=({180:0.5})]A) node[vertex] {$1,a$}
    ([shift=({0:0.5})]A) node[vertex] {$2,a$}
    (2,0) node[cluster] (B) {}
    ([shift=({180:0.5})]B) node[vertex] {$1,b$}
    ([shift=({0:0.5})]B) node[vertex] {$2,b$};
\end{tikzpicture}
\end{center}
Now, for a pair $(Y,G)$ in $X^{1,+}$, the generated subframe $G \subseteq X^1$ is either the cluster on the left, or the cluster on the right, or both of them. Reasoning as above (and taking into consideration condition (iii), which does nothing but avoiding the repetition of clusters that can be \say{contracted one into the other}), we obtain the frame $X^2$ depicted below
\begin{center}
\begin{tikzpicture}[
  vertex/.style = {circle, draw,
  minimum size=26, inner sep=0,
  fill=white},
  cluster/.style = {circle, draw,
  minimum size=22mm, inner sep=0,
  fill=white},
  vertex1/.style = {vertex, fill=red!20!white},
  vertex2/.style = {vertex, fill=orange!30!white},
  vertex3/.style = {vertex, fill=blue!30!white},
  vertex4/.style = {vertex, fill=teal!30!white},
]

  \draw[thick]
    (-2,1.5) node[cluster] (A) {}
    ([shift=({180:0.5})]A) node[vertex] {$1,a$}
    ([shift=({0:0.5})]A) node[vertex] {$2,a$}
    (2,1.5) node[cluster] (B) {}
    ([shift=({180:0.5})]B) node[vertex] {$1,b$}
    ([shift=({0:0.5})]B) node[vertex] {$2,b$}
    (-1,-1.5) node[cluster] (C) {}
    ([shift=({180:0.5})]C) node[vertex] {$1,\delta$}
    ([shift=({0:0.5})]C) node[vertex] {$2,\delta$}
    (1,-1.5) node[vertex] (r) {$\rho,\delta$}
    (-3,-1.5) node[vertex] (a) {$\rho,a$}
    (3,-1.5) node[vertex] (b) {$\rho,b$}
    (C) -- (A)
    (C) -- (B)
    (r) -- (A)
    (r) -- (B)
    (a) -- (A)
    (b) -- (B);
\end{tikzpicture}
\end{center}
\endremark

We will prove that (\ref{eq:prod}) is a diagram of product in $\textsf{S4}\KFr_\lf$.

\proposition\label{prop:belongs}
The cone in \emph{(}\ref{eq:prodn}\emph{)} belongs to $\textup{\textsf{S4}}\KFr_\lf$, for each $n \geq 0$.
\endproposition
\proof
The proof is a simple induction argument on $n$.
\endproof

\remark
For each $n$, the function $X^n \longrightarrow \mathcal{P}(W_0 \times W_1)$ sending $x \in X^n$ to the singleton $\{(p_0^n(x),p_1^n(x))\}$ is a finite irreducible model over the set $W_0 \times W_1$, hence $X$ embeds in the \emph{universal model} for $W_0 \times W_1$ (see the construction of \cite[Section 6]{GdB24}).
\endremark

For $U \in \textsf{S4}\KFr_\lf$, we define $U^n \coloneqq \{u \in U\ \vert\ \depth(u) \leq n\}$ (with the convention that $U^0 = \emptyset$); $U^n \subseteq U$ is a generated subframe, $U^n \subseteq U^{n+1}$ for all $n$ and $U = \cup_n U^n$. Consider a diagram
\[\begin{tikzcd}
	U & {W_0} \\
	& {W_1}
	\arrow["{f_0}", from=1-1, to=1-2]
	\arrow["{f_1}"', from=1-1, to=2-2]
\end{tikzcd}\]
in $\textsf{S4}\KFr_\lf$.
\lemma\label{lem:ext}
With the notations above, given a p-morphism $f^k \colon U^k \longrightarrow X$ such that $p_i f^k = f_i|_{U^k}$ for $i = 0, 1$, there exists a unique p-morphism $f^{k+1} \colon U^{k+1} \longrightarrow X$ such that $f^{k+1}|_{U^k} = f^k$ and $p_i f^{k+1} = f_i|_{U^{k+1}}$ for $i = 0, 1$.
\endlemma
\proof
Assume $f^k \colon U^k \longrightarrow X$ is as required by the hypothesis.
%rev as in the hypothesis being defined.
Consider a cluster $C \subseteq U^{k+1} \setminus U^k$ and define
\begin{align}\label{eq:1}
    Y \coloneqq \{(f_0(u),f_1(u))\ \vert\ u \in C\} \qquad \qquad G \coloneqq f^k(\up{C} \setminus C)\tag{*}
\end{align}
$Y$ is a non-empty finite subset of the cartesian product $W_0 \times W_1$ ($C$ is a nonempty finite set) and $G$ is a finite generated subframe ($\up{C} \setminus C$ is a generated subframe of $U^k$ and the image of a generated subframe via a p-morphism is a generated subframe; $\up{C} \setminus C$ is finite, $U$ being locally-finite and $C$ being finite) of $X^n$ (where $n \coloneqq \depth(G)$), comprising at least a point of depth $n$, if $n > 0$. Moreover, given $y \in Y$, say $y = (f_0(u),f_1(u))$ for $u \in C$,
\begin{align*}
    \pi_i(Y) \cup p_i^n(G) &= \{f_i(u')\}_{u' \in C} \cup f_i(\up{C} \setminus C) = \\
    &= f_i(\up{C}) = f_i(\up{u}) = \up{f_i(u)} = \up{\pi_i(y)}
\end{align*}
for all $i = 0, 1$, $f_i$ being a p-morphism and $U$ being in $\textsf{S4}\KFr_\lf$. This proves conditions (i) and (ii) in the construction of $X^{n,+}$ above for the pair $(Y,G)$.

Depending on whether condition (iii) is satisfied or not, we define $f^{k+1}$ on $C$ as follows.
\begin{enumerate}
    \item[a)] If $(Y,G) \in X^{n,+}$, we define, for $u \in C$, $f^{k+1}(u)$ to be the point $((f_0(u),f_1(u)),Y,G)$ in the cluster $\{(y,Y,G)\ \vert\ y \in Y\}$; with such definition,
    \begin{align*}
        p_i f^{k+1}(u) = p_i^{n+1}(((f_0(u),f_1(u)),Y,G)) = \pi_i((f_0(u),f_1(u))) = f_i(u)
    \end{align*}
    for all $u \in C$.
    \item[b)] If $(Y,G) \notin X^{n,+}$, i.e. $\{(f_0(u),f_1(u))\ \vert\ u \in C\} \subseteq \{(p_0^n(x),p_1^n(x))\ \vert\ x \in X^n \text{ s.t.\ } \up{x} = G\}$ we define, for $u \in C$, $f^{k+1}(u)$ to be the unique $x \in X^n$ such that $\up{x} = G$ and $(f_0(u),f_1(u)) = (p_0^n(x),p_1^n(x))$ (if $D \subseteq X^n$ is a cluster of depth $n$, then the map $D \ni x \mapsto (p_0^n(x),p_1^n(x))$ is injective); with such definition,
    \begin{align*}
        p_i f^{k+1}(u) = p_i^n f^{k+1}(u) = f_i(u)
    \end{align*}
    for all $u \in C$.
\end{enumerate}
If we repeat the same scheme for each cluster $C \subseteq U^{k+1} \setminus U^k$, and we set $f^{k+1}|_{U^k} = f^k$, we obtain a function $f^{k+1} \colon U^{k+1} \longrightarrow X$, which extends $f^k$ and such that $p_i f^{k+1} = f_i|_{U^{k+1}}$ (the set $U^{k+1} \setminus U^k$ can be partitioned into clusters).

The function $f^{k+1}$ just defined is a p-morphism. Consider $u \in U^{k+1}$. If $u \in U^k$, then $\up{u} \subseteq U^k$, hence
\begin{align*}
    f^{k+1}(\up{u}) = f^k(\up{u}) = \up{f^k(u)} = \up{f^{k+1}(u)}
\end{align*}
since $f^k$ is a p-morphism and $f^{k+1}|_{U^k} = f^k$. If $u \in U^{k+1} \setminus U^k$, then there exists a unique cluster $C \subseteq U^{k+1} \setminus U^k$ such that $u \in C$; taking $Y$ and $G$ as in (\ref{eq:1}), we have that
\begin{enumerate}
    \item[a)] if $(Y,G) \in X^{n,+}$, then $\up{f^{k+1}(u)} = \{(y,Y,G)\ \vert\ y \in Y\} \cup G$, hence
        \begin{align*}
            f^{k+1}(\up{u}) &= f^{k+1}(\up{C}) = f^{k+1}(C) \cup f^k(\up{C} \setminus C) = \\
            &= \{(y,Y,G)\ \vert\ y \in Y\} \cup G = \up{f^{k+1}(u)}
        \end{align*}
    \item[b)] if $(Y,G) \notin X^{n,+}$, i.e.\ $Y \subseteq \{(p_0^n(x),p_1^n(x))\ \vert\ x \in X^n \text{ s.t.\ } \up{x} = G\}$, then $\up{f^{k+1}(u)} = G$, hence
        \begin{align*}
            f^{k+1}(\up{u}) &= f^{k+1}(\up{C}) = f^{k+1}(C) \cup f^k(\up{C} \setminus C) = \\
            &= f^{k+1}(C) \cup G = G = \up{f^{k+1}(u)}
        \end{align*}
\end{enumerate}

Now, consider another extension $g^{k+1}$ of the p-morphism $f^k$ as in our claim. We have to verify that, given a cluster $C \subseteq U^{k+1} \setminus U^k$ as above, $f^{k+1}(u) = g^{k+1}(u)$ for all $u \in C$; observe that, since the two extensions satisfy $p_i f^{k+1}(u) = f_i(u) = p_i g^{k+1}(u)$, for $i = 0, 1$, for all $u \in C$, it is sufficient to prove that $f^{k+1}(C)$ and $g^{k+1}(C)$ are contained in the same cluster (if $D \subseteq X$ is a cluster, then $D \ni x \mapsto (p_0(x),p_1(x))$ is injective, by construction of $X$). Recall that we defined
\begin{align*}
    Y = \{(f_0(u),f_1(u))\ \vert\ u \in C\} \qquad \qquad G = f^k(\up{C} \setminus C)
\end{align*}
and $n = \depth(G)$. As previously observed, the pair $(Y,G)$ satisfies conditions (i) and (ii) in the definition of $X^{n,+}$, hence $(Y,G) \in X^{n,+}$ if and only if $(Y,G)$ satisfies condition (iii). Since $p_i g^{k+1} = f_i|_{U^{k+1}}$, for $i = 0, 1$, and $g^{k+1}|_{U^k} = f^k$, we have
\begin{align*}
    Y = \{(p_0(x),p_1(x))\ \vert\ x \in g^{k+1}(C)\} \qquad \qquad G = g^{k+1}(\up{C} \setminus C)
\end{align*}
In particular, $\depth(g^{k+1}(C))$ can be $n$ (if $g^{k+1}(C) \subseteq G$), or $n+1$ (otherwise); in the latter case, $g^{k+1}(C)$ must coincide with the cluster it is contained in.
\begin{enumerate}
    \item[a)] If $(Y,G) \in X^{n,+}$, then $g^{k+1}(C) \nsubseteq G$, otherwise we would have that $\up{g^{k+1}(C)} = G$ and, as a consequence,
    \begin{align*}
    Y &= \{(p_0(x),p_1(x))\ \vert\ x \in g^{k+1}(C)\} \subseteq \\
    &\subseteq \{(p_0^n(x),p_1^n(x))\ \vert\ x \in X^n \text{ s.t. } \up{x} = G\}
    \end{align*}
    against $(Y,G) \in X^{n,+}$. This means that $g^{k+1}(C)$ is a cluster in $X$ of depth $n+1$, say of the form $\{(y',Y',G')\ \vert\ y' \in Y'\}$ for some $(Y',G') \in X^{n,+}$. We have that
    \begin{align*}
    G' = \up{g^{k+1}(C)} \setminus g^{k+1}(C) = g^{k+1}(\up{C} \setminus C) = G
    \end{align*}
    by definition of the binary relation in $X^{n+1}$ and by $g^{k+1}(C) \nsubseteq G$, which is equivalent to $g^{k+1}(C) \cap G = \emptyset$, and
    \begin{align*}
    Y' = \{(\pi_0(y'),\pi_1(y'))\ \vert\ y' \in Y'\} = \{(p_0(x),p_1(x))\ \vert\ x \in g^{k+1}(C)\} = Y
    \end{align*}
    \item[b)] If $(Y,G) \notin X^{n,+}$, then $g^{k+1}(C) \subseteq G$, otherwise $g^{k+1}(C)$ would be a cluster in $X$ of depth $n+1$, say $g^{k+1}(C) = \{(y',Y',G')\ \vert\ y' \in Y'\}$ for some $(Y',G') \in X^{n,+}$, and, reasoning as above, this would imply $(Y,G) = (Y',G') \in X^{n,+}$. As a consequence,
    $$\up{g^{k+1}(C)} = G = \up{f^{k+1}(C)}.$$
\end{enumerate}
In both cases, $g^{k+1}(C)$ and $f^{k+1}(C)$ are contained in the same cluster.
\endproof

\proposition\label{prop:prod}
With the notation above, $(p_i \colon X \longrightarrow W_i\ \vert\ i = 0, 1)$ is the product of $W_0$ and $W_1$ in $\textup{\textsf{S4}}\KFr_\lf$.
\endproposition
\proof
The cone $(p_i \colon X \longrightarrow W_i\ \vert\ i = 0, 1)$ is in $\textsf{S4}\KFr_\lf$, since $X$ is the colimit in $\KFr$ of the diagram of the inclusions $X^n \longrightarrow X^m$ (which lives in $\textsf{S4}\KFr_\lf$, by Proposition~\ref{prop:belongs}) and the embedding $\textsf{S4}\KFr_\lf \longrightarrow \KFr$ creates colimits (by Proposition~\ref{prop:colim}).% The morphisms $p_0$ and $p_1$ are induced by the universal property of the colimit.

Consider another cone $(f_i \colon U \longrightarrow W_i\ \vert\ i = 0, 1)$ in $\textsf{S4}\KFr_\lf$ and let $f^0 \colon U^0 = \emptyset \longrightarrow X$ be the empty function. Using Lemma~\ref{lem:ext}, we can inductively define a p-morphism $f^k \colon U^k \longrightarrow X$ such that $p_i f^k = f_i|_{U^k}$ for $i = 0, 1$, and $f^{k+1}|_{U^k} = f^k$ for $k \geq 0$. Since each $U^k$ is a generated subframe of $U$ and $U = \cup_n U^n$, the $f^k$'s can be glued to a p-morphism $f \colon U \longrightarrow X$ such that $p_i f = f_i$ for $i = 0, 1$. Now, given another p-morphism $g \colon U \longrightarrow X$ such that $p_i g = f_i$ for $i = 0, 1$, and setting $\bar{k} \coloneqq \text{min}\{k \geq 0\ \vert\ f|_{U^k} \neq g|_{U^k}\}$ ($\bar{k} > 0$, since both $f|_{U^0}$ and $g|_{U^0}$ must be the empty function), applying the uniqueness part of Lemma~\ref{lem:ext} to the extension of the p-morphism $f|_{U^{\bar{k}-1}} = g|_{U^{\bar{k}-1}}$ (they are equal by the minimality of $\bar{k}$), we conclude that $f|_{U^{\bar{k}}} = g|_{U^{\bar{k}}}$ (against the definition of $\bar{k}$).
\endproof

To generalize the description of the coproduct of profinite algebras for any $L \in \NExt \textsf{S4}$ with the finite model property, we have to add some constraints in the construction of $X^{n,+}$. Consider the construction above for $W_0$ and $W_1$ in $L\KFr_\lf$.

\proposition\label{prop:prodc}
With the notations above, the subset $X_L \coloneqq \{x \in X\ \vert\ x^* \in L\KFr_\lf\} \subseteq X$ is a generated subframe of $X$ and, if we set $(p_i)_L \coloneqq p_i|_{X_L}$, $((p_i)_L \colon X_L \longrightarrow W_i\ \vert\ i = 0, 1)$ is the product of $W_0$ and $W_1$ in $L\KFr_\lf$.
\endproposition
\proof
It is sufficient to observe that, given a cone $(f_i \colon U \longrightarrow W_i\ \vert\ i = 0, 1)$ in $L\KFr_\lf$, the unique p-morphism $f \colon U \longrightarrow X$ induced by the universal property of the product in $\textsf{S4}\KFr_\lf$ (see Proposition~\ref{prop:prod}) factors through $X_L$. For, given $u \in U$, $f(u)^* = f(u^*)$ belongs to $L\KFr_\lf$, being a p-morphic image of $u^*$,\footnote{If $W \longrightarrow V$ is a surjective p-morphism, then the logic of $W$ is contained in the logic of $V$ (see \cite{misha}).} which is a generated subframe of $U \in L\KFr_\lf$.
\endproof

\begin{comment}
\remark\label{rem:pos}
If we take $W_0$ and $W_1$ in $\textit{Grz}\KFr_\lf = \{\text{locally finite posets}\}$ and we compute their product in $S4\KFr_\lf = \{\text{locally finite preorders}\}$, then the result belongs to $\textit{Grz}\KFr_\lf$, giving the product of $W_0$ and $W_1$ in $\textit{Grz}\KFr_\lf$. In fact, if $(Y,G) \in X^{n,+}$ and $\lvert Y \rvert \geq 2$ (say $(w_0,w_1), (w_0',w_1') \in Y$, with $w_0 \neq w_0'$), then we would have
\begin{align*}
     \up{w_0} = \up{\pi_0(w_0,w_1)} = \pi_0(Y) \cup p_0^n(G) = \up{\pi_0(w_0',w_1')} = \up{w_0'}
\end{align*}
hence $w_0 = w_0'$, which is a contradiction.
\endremark
\end{comment}

\section{Coregularity}\label{sec:AP}
As we said in the introduction, in \cite{infinitary} we developed infinitary sequent calculi for infinitary propositional modal languages, by adding specific modal rules to the Maehara-Takeuti infinitary extension of \textbf{LK} (see \cite{Tak}). We then recovered profinite $L$-algebras as Lindenbaum algebras for such calculi, connecting syntactic properties to categorical ones. In particular, we proved that the infinitary Craig's Interpolation Theorem for the global consequence relation corresponds to the amalgamability of $\Pro L\MA_\fin$. If all epimorphisms in $\Pro L\MA_\fin$ are regular (for example, whenever $L \in \NExt \textsf{K4}$), which happens exactly when an infinitary version of the standard Beth's Definability Property holds, we can say more: $\Pro L\MA_\fin$ satisfies the amalgamation property if and only if $L\KFr_\lf \simeq (\Pro L\MA_\fin)^\op$ is a regular category.

The main goal of this section is to characterize those normal modal logics $L \in \NExt \textsf{S4}$ with the finite model property for which the category $L\KFr_\lf$ is regular. We then compare this characterization with Maksimova's classification of the amalgamability of $L\MA$, which corresponds to $L$ satisfying the Craig's Interpolation Theorem for necessary propositions \cite{Maksimova1979InterpolationTI,Maksimova1980InterpolationTI}.

The section is organized as follows. We first establish the equivalence between regularity of $L\KFr_\lf$ and amalgamability of $\Pro L\MA_\fin$ and $L\MA_\fin$, under the assumption that $L \in \NExt \textsf{K4}$ (the proof of this result can also be found in \cite{infinitary}). We then move to the combinatorics of finite Kripke frames and derive the necessary restrictions on $L$. Finally, we prove that all the resulting logics satisfy the required property, obtaining the desired classification.

\subsection{Coregularity and amalgamability}
Let us recall the notion of regularity (see \cite{MR1953060}), in one of its versions that uses \emph{extremal} epimorphisms.\footnote{An epimorphism is said to be \emph{extremal} if it does not factorize through any proper subobject of its codomain.}
\definition\label{defn:regular}
A category $\mathbf{C}$ is \emph{regular} if
\begin{enumerate}
    \item it has finite limits;
    \item it has (extremal epi, mono)-factorization;
    \item extremal epimorphisms are pullback-stable.
\end{enumerate}
\enddefinition

The \emph{amalgamation property} is usually formulated for categories of algebras, but it can be formulated in purely categorical terms.
\definition
We say that a category $\mathbf{C}$ satisfies the \emph{amalgamation property} \emph{(AP)} if, for any monomorphisms $k \colon A \longrightarrow B$, $g \colon A \longrightarrow C$, there are monomorphisms $h \colon C \longrightarrow D$, $f \colon B \longrightarrow D$ such that $f k = h g$, i.e.\ such that the square
\[\begin{tikzcd}
	A & B \\
	C & D
	\arrow["k", from=1-1, to=1-2]
	\arrow["g"', from=1-1, to=2-1]
	\arrow["f", from=1-2, to=2-2]
	\arrow["h"', from=2-1, to=2-2]
\end{tikzcd}\]
commutes. The pair of monomorphisms $h$ and $f$ is called \emph{amalgam} for the pair $k$ and $g$. The dual of the amalgam, namely a commutative square of epimorphisms, is called \emph{coamalgam}.
\enddefinition

Let $L$ be a normal modal logic with the finite model property. We want to describe (extremal) epimorphisms and monomorphisms in $L\KFr_\lf \simeq (\Pro L\MA_\fin)^\op$. Consider a p-morphism $f \colon W \longrightarrow V$ in $L\KFr_\lf$. Define
$$U \coloneqq \{(v,i) \in V \times \{0,1\}\ \vert\ \text{if } i = 0, \text{ then } v \notin \Image f\}$$
with inclusions $\iota_i \colon V \longrightarrow U$ such that $\iota_0(v) = (v,0)$ if $v \notin \Image f$, $\iota_0(v) = (v,1)$ if $v \in \Image f$, and $\iota_1(v) = (v,1)$. Observe that $U$ is obtained by gluing two copies of $V$ along $\Image f$ and that $\iota_0$ and $\iota_1$ are the corresponding inclusions: this means that $\iota_0$ and $\iota_1$ is the cokernel pair of $f$ in $\Set$. To obtain the cokernel pair in $\KFr$, hence in $L\KFr_\lf$ (see Proposition~\ref{prop:colim}), it is necessary to endow $U$ with the following binary relation:
$$(v,i) \prec (v',i') ~~~~~\iff~~~~~ \begin{cases}
    v \prec v' \text{ and } i = i',\\
    \text{or}\\
    v \prec v' \text{ and } i<i' \text{ and } v' \in \Image f.
\end{cases}$$
It is straightforward to verify that the relation above makes $\iota_0$ and $\iota_1$ into p-morphisms; moreover, given any other pair of parallel p-morphisms coequalizing $f$, the function induced by the universal property of $\iota_0$ and $\iota_1$ as cokernel pair of $f$ in $\Set$ is a p-morphism. We conclude that $\iota_0$ and $\iota_1$ is the cokernel pair of $f$ in $L\KFr_\lf$.
%rev With such definitions, $\iota_0$ and $\iota_1$ are p-morphisms and $U \in L\KFr_\lf$ ($U$ satisfies $L$, being a p-morphic image of the disjoint union $V + V$, see \cite{misha}). Moreover, it is straightforward to verify that the pair $\iota_0, \iota_1$ is the cokernel pair of $f$ in $L\KFr_\lf$.
The equalizer in $L\KFr_\lf$ of $\iota_0$ and $\iota_1$ can be easily identified with the image of $f$, since $\iota_0(v) = \iota_1(v)$ if and only if $v \in \Image f$. This construction gives us the so-called\emph{ coregular factorization} in $L\KFr_\lf$
\[\begin{tikzcd}[column sep=small]
	W && V && U \\
	& {\Image f}
	\arrow["f", from=1-1, to=1-3]
	\arrow["e"', from=1-1, to=2-2]
	\arrow["{\iota_1}"', shift right, from=1-3, to=1-5]
	\arrow["{\iota_0}", shift left, from=1-3, to=1-5]
	\arrow["m"', from=2-2, to=1-3]
\end{tikzcd}\]
where $m$ is the inclusion of $\Image f$ in $V$ and $e$ sends $w \in W$ to $f(w) \in \Image f$ (it is the unique morphism in $L\KFr_\lf$ induced by the universal property of the equalizer). Before stating the next results, we make this simple observation: the isomorphisms in $L\KFr_\lf$ are exactly the bijective p-morphisms.
\proposition\label{prop:injsurj}
Let $L$ be a normal modal logic with the finite model property. A morphism $f$ in $L\KFr_\lf$ is
\begin{enumerate}
    \item an epimorphism if and only if it is surjective;
    \item a regular monomorphism if and only if it is injective.
\end{enumerate}
%Moreover, if $L$ extends $K4$, then all monomorphisms are regular in $L\KFr_\lf$.
\endproposition
\proof
For 1. observe that $f$ is an epimorphism if and only if $\iota_0 = \iota_1$, if and only if $m$ is an isomorphism, if and only if $f$ is surjective.

For 2. recall that an arrow is a regular monomorphism iff it is the equalizer of its own cokernel pair; thus $f$ is a regular monomorphism if and only if $e$ is an isomorphism, if and only if $f$ is injective.
\endproof

If $L$ extends \textsf{K4}, we can say more.
\proposition\label{prop:injsurj2}
Let $L$ be a normal modal logic with the finite model property and assume that $L \in \NExt \textup{\textsf{K4}}$. A morphism $f$ in $L\KFr_\lf$ is
\begin{enumerate}
    \item an extremal epimorphism if and only if it is surjective;
    \item a monomorphism if and only if it is injective.
\end{enumerate}
\endproposition
\proof
The right-to-left implication in 2. follows from Proposition~\ref{prop:injsurj}. For the other direction, let $f \colon W \longrightarrow V$ be a monomorphism in $L\KFr_\lf$. If $f$ is not injective, then we have some $w \in W$ for which there exists $w' \in W$ with $w' \neq w$ and $f(w') = f(w)$; consider $w$ having minimal $|w^*|$ with this property.

Consider first the case where $w' \in w^*$; by the minimality condition with respect to the property above, it must be $w \in w'^*$ (otherwise $|w'^*| < |w^*|$), hence $w^*=w'^*$. Thanks to transitivity (here we are using the assumption that $L$ contains \textsf{K4}), the function $\sigma_{w,w'} \colon w^* \longrightarrow w^*$ swapping $w$ and $w'$ (and fixing all the other points) is a p-morphism in $L\KFr_\lf$ ($w^*$ is a generated subframe of $W$). By $f(w)=f(w')$, the following diagram in $L\KFr_\lf$
\[\begin{tikzcd}
	{w^*} & W \\
	W & V
	\arrow["\iota", from=1-1, to=1-2]
	\arrow["{\iota \sigma_{w,w'}}"', from=1-1, to=2-1]
	\arrow["f", from=1-2, to=2-2]
	\arrow["f"', from=2-1, to=2-2]
\end{tikzcd}\]
commutes ($\iota$ is the inclusion $w^* \subseteq W$). $f$ being a monomorphism, we have that $\iota \sigma_{w,w'} = \iota$; in particular, $w = w'$, which is a contradiction.

The contradiction reached for the case $w' \in w^*$ implies that the p-morphism $w^* \longrightarrow f(w^*) = f(w)^*$ in $L\KFr_\lf$, restriction of $f$, is injective; being surjective, it is an isomorphism. We can then consider its inverse $h \colon f(w^*) \longrightarrow w^*$. The morphism $h$ is now used to prove that the case $w' \notin w^*$ gives a contradiction, too.

Assume that $w' \notin w^*$. Since $f(w) = f(w')$, we have $f(w^*) = f(w)^* = f(w')^* = f(w'^*)$, so we can consider the p-morphism $g \colon w'^* \longrightarrow w^*$ in $L\KFr_\lf$, given by the composition of the restriction $w'^* \longrightarrow f(w'^*)$ of $f$ with $h \colon f(w^*) \longrightarrow w^*$. By definition, the following diagram in $L\KFr_\lf$
\[\begin{tikzcd}
	{w'^*} & W \\
	W & V
	\arrow["{\iota'}", from=1-1, to=1-2]
	\arrow["{\iota g}"', from=1-1, to=2-1]
	\arrow["f", from=1-2, to=2-2]
	\arrow["f"', from=2-1, to=2-2]
\end{tikzcd}\]
commutes ($\iota'$ is the inclusion $w'^* \subseteq W$). $f$ being a monomorphism, we have that $\iota g = \iota'$; in particular, $w = w'$, which is a contradiction.

The left-to-right implication in 1. follows from Proposition~\ref{prop:injsurj}. For the other direction, assume that $f$ is surjective and consider a factorization $f = m e$, with $m$ mono. The p-morphism $m$ is then surjective and injective, by 2., hence it is an isomorphism.
\endproof

We are now ready to connect regularity and amalgamability, under the assumption that $L \in \NExt \textsf{K4}$. In view of Propositions~\ref{prop:injsurj} and~\ref{prop:injsurj2}, we refer to epimorphisms, regular epimorphisms and surjective p-morphisms simply as ‘epimorphisms’; similarly, we refer to monomorphisms, regular monomorphisms and injective p-morphisms simply as ‘monomorphisms’.
\theorem\label{teo:regular}
Let $L$ be a normal modal logic with the finite model property and assume that $L \in \NExt \textup{\textsf{K4}}$. Then $L\KFr_\lf$ is regular if and only if $\Pro L\MA_\fin$ satisfies the \emph{(AP)}.
\endtheorem
\proof
The category $L\KFr_\lf$ enjoys properties (i) and (ii) in Definition~\ref{defn:regular} of regular category (by Corollary~\ref{cor:lim} and by the existence of the coregular factorization). Moreover, the pullback of an epimorphism along any monomorphism is an epimorphism: it is given by the inverse image. By applying the pullback lemma to the image factorization we can then reduce regularity to the following: epimorphisms are stable under pullback along epimorphisms. The claim now follows from a general fact in category theory, recalling that $L\KFr_\lf$ is dual to $\Pro L\MA_\fin$: a category with pushouts satisfies the (AP) if and only if monomorphisms are stable under pushout along monomorphisms.
\endproof

Recall that $\Pro L\MA_\fin$ contains the class $L\MA_\fin$ of finite $L$-algebras.
\proposition\label{prop:APfin}
Let $L$ be a normal modal logic with the finite model property. Then $\Pro L\MA_\fin$ satisfies the \emph{(AP)} if and only if $L\MA_\fin$ satisfies the \emph{(AP)}.
\endproposition
\proof
We consider the dual (AP) for $L\KFr_\lf \simeq (\Pro L\MA_\fin)^\op$ and $L\KFr_\fin \simeq (L\MA_\fin)^\op$ (in $L\KFr_\fin$, as well as in $L\KFr_\lf$, the epimorphisms are the surjective p-morphisms: the proof works as in Proposition~\ref{prop:injsurj}). Consider a diagram
\[\begin{tikzcd}
	& {U} \\
	{W} & V
	\arrow["{g}", from=1-2, to=2-2]
	\arrow["{f}"', from=2-1, to=2-2]
\end{tikzcd}\]
of Kripke frames validating $L$, where both $f$ and $g$ are surjective p-morphisms.

One direction is straightforward. If the diagram above lives in $L\KFr_\fin$ and we have a coamalgam in $L\KFr_\lf$, say $h \colon A \longrightarrow W$ and $k \colon A \longrightarrow U$, then we also get a coamalgam in $L\KFr_\fin$ for the same diagram. For, restrict $h$ and $k$ to $A' \coloneqq (\bigcup a_w^*) \cup (\bigcup a_u^*) \subseteq A$, where, for each $w \in W$, a point $a_w \in A$ such that $h(a_w) = w$ has been fixed (and the same for each $u \in U$). Then $A' \in L\KFr_\fin$: it is finite because $W + U$ is finite and $A$ is locally finite; it satisfies $L$ because $A$ does so and $A'$ is a generated subframe (taking generated subframes is logic-preserving \cite{misha}).

For the other direction, assume that the diagram above lives in $L\KFr_\lf$. Let %an element of the $\DGrph$-pullback $P_{f,g}$, i.e.\
$(w,u)$ be a pair in the cartesian product $W \times U$ and assume that $f(w) = g(u)$; we can consider the generated subframe
$$V(w,u) \coloneqq f(w^*) = f(w)^* = g(u)^* = g(u^*)$$
of $V$ and the restricted diagram
\[\begin{tikzcd}
	& {u^*} \\
	{w^*} & {V(w,u)}
	\arrow["{g_{u}}", from=1-2, to=2-2]
	\arrow["{f_{w}}"', from=2-1, to=2-2]
\end{tikzcd}\]
in $L\KFr_\fin$ involving surjective p-morphisms ($f_{w}$ is the restriction of $f$ and $g_{u}$ is the restriction of $g$). By hypothesis, we can find some coamalgam
\[\begin{tikzcd}
	{A(w,u)} & {u^*} \\
	{w^*} & {V(w,u)}
	\arrow["{k_{w,u}}", from=1-1, to=1-2]
	\arrow["{h_{w,u}}"', from=1-1, to=2-1]
	\arrow["{g_{u}}", from=1-2, to=2-2]
	\arrow["{f_{w}}"', from=2-1, to=2-2]
\end{tikzcd}\]
in $L\KFr_\fin$. Take the disjoint union
$$A \coloneqq \{(w,u,a)\ \vert\ w \in W,\ u \in U \text{ s.t. } f(w) = g(u) \text{ and } a \in A(w,u)\},$$ %$\coprod A(w,u)$, varying $(w,u)$ among those pairs such that $f(w) = g(u)$.
with
$$(w,u,a) \prec (w',u',a') ~~~~~\iff~~~~~ w = w',\ u = u' \text{ and } a \prec a' \text{ in } A(w,u).$$
The induced $h \colon A \longrightarrow W$ and $k \colon A \longrightarrow U$ are p-morphisms in $L\KFr_\lf$ (taking the disjoint union of Kripke frames is logic-preserving \cite{misha}). Moreover, $h$ and $k$ are surjective. For, let $w \in W$: by the surjectivity of $g$, we can find $u \in U$ such that $f(w) = g(u)$ and then, by the surjectivity of $h_{w,u}$, there exists $a \in A(w,u)$ such that $h(w,u,a) = h_{w,u}(a) = w$. The surjectivity of $k$ follows by symmetry. Since $f h = g k$, we conclude that the pair $h$ and $k$ is a coamalgam in $L\KFr_\lf$ for the given diagram.
\endproof

\begin{comment}
Let $L \in \NExt \textsf{S4}$. To $L$ we can associate a \emph{superintuitionistic logic} (\emph{si-logic} for short), i.e.\ an extension of the propositional intuitionistic logic \textsf{Ipc} (we denote by $\Ext \textsf{Ipc}$ the family of si-logics): such si-logic is denoted by $\rho L$ and it is called the
\emph{superintuitionistic fragment} of $L$ (see \textcolor{red}{reference}); $L$ is said to be a \emph{modal companion} of $\rho L$. Vice versa, given $I \in \Ext \textsf{Ipc}$, then there exists $\tau I \in \NExt \textsf{S4}$ such that, for all $L \in \NExt \textsf{S4}$,
$$\rho L = I ~~~~~\iff~~~~~ \tau I \subseteq L \subseteq \tau I + \textsf{Grz}.$$

Given $I \in \Ext \textsf{Ipc}$, we denote by $I\HA$ the category of the algebraic models of $I$ and by $I\HA_\fin$ the full subcategory of the finite members of $I\HA$.

\proposition\label{prop:sifragment}
Let $L \in \NExt \textup{\textsf{S4}}$. If $L\MA$ has the \emph{(AP)}, then $\rho L \HA$ has the \emph{(AP)}. If $L\MA_\fin$ has the \emph{(AP)}, then $\rho L \HA_\fin$ has the \emph{(AP)}.
\endproposition
\proof
It follows from \cite[Proposition 2.1]{Maksimova1979InterpolationTI} and \cite[\S 3, Theorem 3.1]{MR163850}.
\endproof
\end{comment}

The (AP) can be encountered in many cases. Given p-morphisms $f \colon W \longrightarrow V$ and $g \colon U \longrightarrow V$, we denote by $P_{f,g}$ the pullback in $\DGrph$\footnote{Recall that $\DGrph$ is the category of Kripke frames and stable functions, i.e.\ satisfying $w \prec w' \implies f(w) \prec f(w')$.} of $f$ and $g$, i.e.\ the Kripke frame having
$$\{(w,u) \in W \times U\ \vert\ f(w) = g(u)\}$$
as underlying set and equipped with the product relation
$$(w,u) \prec (w',u') ~~~~~\iff~~~~~ w \prec w' \text{ and } u \prec u'.$$
The two projections $\pi_W \colon P_{f,g} \longrightarrow W$ and $\pi_U \colon P_{f,g} \longrightarrow U$ are p-morphisms. Moreover, $\pi_W$ is surjective if $g$ is surjective (symmetrically, $\pi_U$ is surjective if $f$ is surjective). However, even if $W$ and $U$ validate $L$, it is not guaranteed that $P_{f,g}$ does so.

\remark\label{rem:Horn}
Assume that $L\KFr_\fin$ can be characterized via \emph{universal Horn clauses} in the first-order language containing the binary predicate $\prec$ and equality.\footnote{A universal Horn clause is a disjunction of \textit{literals} --- formulas of the kind $R(x_1,\dots,x_n)$ or $\neg R(x_1,\dots,x_n)$, for some first-order relational symbol $R$ --- in which at most one of them is unnegated.} It is straightforward to verify that, under these assumptions, if $f \colon W \longrightarrow V$ and $g \colon U \longrightarrow V$ are surjective p-morphisms in $L\KFr_\fin$, then the Kripke frame $P_{f.g}$ validates $L$; as a consequence, $\pi_W$ and $\pi_U$ are a coamalgam for $f$ and $g$ in $L\KFr_\fin$. Dualizing, we obtain that $L\MA_\fin$ satisfies the (AP). This observation can be applied to the logics in the Table in Section~\ref{sec:profinite}: if $L$ is one among \textsf{K}, \textsf{K4}, \textsf{S4}, \textsf{S4.2}, \textsf{Grz}, \textsf{S5} and \textsf{GL}, then $L\MA_\fin$ and $\Pro L\MA_\fin$ satisfy the (AP), and $L\KFr_\lf$ is a regular category.
\endremark

We conclude this subsection by comparing the amalgamability of $L\MA_\fin$ with the amalgamability of the full variety $L\MA$.
\proposition\label{prop:APalg}
Let $L$ be a normal modal logic with the finite model property and assume that $L \in \NExt \textup{\textsf{K4}}$. If $L\MA$ satisfies the \emph{(AP)}, then $L\MA_\fin$ satisfies the \emph{(AP)}.
\endproposition
\proof
Apply \cite[Propositions 2.16 and 2.20]{GZ}, observing that the equational theory of $L\MA$ has \textit{equationally definable operations} (see \cite[Chapter 2.5]{GZ}).
\endproof

\subsection{Necessary conditions}
%We showed that, if $L \in \NExt \textsf{K4}$, then the category $L\KFr_\lf \simeq (\Pro L\MA_\fin)^\op$ is regular if and only if $\Pro L\MA_\fin$ satisfies the (AP), if and only if $L\MA_\fin$ satisfies the (AP).
In \cite{Maksimova1979InterpolationTI}, Maksimova provided a finite list containing all the logics $L \in \NExt \textsf{S4}$ such that $L\MA$ satisfies the (AP). In this subsection, we prove that the same list exhausts all cases above \textsf{S4} for which $L\MA_\fin$ satisfies the (AP), under the additional assumption that $L$ has the finite model property.\footnote{This assumption is non-restrictive, since it is satisfied by all logics occurring in Maksimova's classification.} To do so, we use the finite duality $L\MA_\fin \simeq (L\KFr_\fin)^\op$ and combinatorial methods on finite Kripke frames, following \cite{GZ}. The logics of the resulting classification will be presented by means of the classes of their finite Kripke frames. Indeed, the finite model property ensures that a logic $L$ is uniquely determined by the class $L\KFr_\fin$ of finite Kripke frames validating $L$. Moreover, a class of finite transitive Kripke frames is modally definable (meaning that it coincides with $L\KFr_\fin$, for some $L \in \NExt \textsf{K4}$) if and only if it is closed under generated subframes, p-morphic images and finite disjoint unions (see \cite[Theorem 3.21]{MR1837791}). 

The following invariants serve to define classes of (locally) finite Kripke frames satisfying the closure properties mentioned above. Recall that, for a locally finite preorder $W$ and $w \in W$, the \emph{depth} $\depth(w)$ is defined as the maximum length of sequences $w = w_1 \prec w_2 \prec \cdots \prec w_n$, with $w_i \neq w_{i+1}$ and $w_{i+1} \nprec w_i$. Moreover, we define the \emph{width} $\width(w)$ as the maximum number of mutually incomparable points in $w^*$. Finally, let $\numberext(w)$ denote the number of external clusters in $w^*$. Then, for each $\alpha \in \{\depth, \width, \numberext\}$, we define $\alpha(W)$ as the maximum of $\{\alpha(w)\ \vert\ w \in W\}$, with the convention that $\alpha(W) = \omega$ if the maximum does not exist and $\alpha(W) = 0$ if $W = \emptyset$.

The following result provides a preliminary classification in a restricted setting. Recall that $\textsf{Grz}\KFr_\lf$ is the class of locally finite partial orders.
\proposition\label{prop:regnecpos}
Let $L$ be a normal modal logic with the finite model property and assume that $L \in \NExt \textup{\textsf{Grz}}$. If $L\MA_\fin$ satisfies the \emph{(AP)}, then $L = \textup{\textsf{Grz.3}}$, or $L\KFr_\fin$ is one of the following:
\begin{enumerate}
    \item[(I)] $\{W \in \textup{\textsf{Grz}}\KFr_\fin\ \vert\ \depth(W) \leq l,\ \width(W) \leq k\}$, for some $l \in \{0,1,2\}$ and some $k \in \{1,2,\omega\}$;
    \item[(II)] $\{W \in \textup{\textsf{Grz}}\KFr_\fin\ \vert\ \numberext(W) \leq s\}$, for some $s \in \{1,\omega\}$.
\end{enumerate}
\endproposition
\proof
If $L$ is as in the hypothesis, then $L\MA_\fin$ is equivalent to the category of finite algebras of the variety of Heyting algebras corresponding to $L$ via Block-Esakia, \cite{MR3967703}: the transformations witnessing Block-Esakia Theorem (see \cite{MR1171642} for a survey) define a categorical equivalence between finite Heyting algebras and finite \textsf{Grz}-algebras and such an equivalence restricts to subvarieties. The classes of finite partial orders listed above are precisely those corresponding to the amalgamable varieties of Heyting algebras; the claim then follows from the fact that, for varieties of Heyting algebras, finite amalgamability (i.e.\ finite algebras being amalgamable) implies the amalgamability of the whole variety, \cite[Theorem 2]{Maksimova1977CraigsTI}.
\endproof

The next result will be useful for the remainder of this section. Given a preorder $W$, we denote by $\Refl W$ its \emph{posetal reflection}, i.e.\ the set of clusters $x \subseteq W$, ordered by
$$x \prec x' ~~~~~\iff~~~~~ \text{there are } w \in x \text{ and } w' \in x' \text{ s.t. } w \prec w' \text{ in } W.$$
In other words, $\Refl W$ is the quotient of $W$ via the equivalence relation
$$w \sim w' ~~~~~\iff~~~~~ w \prec w' \text{ and } w' \prec w.$$
It is straightforward to prove that $\Refl W$ is a partial order. We have an obvious projection $\eta_W \colon W \longrightarrow \Refl W$, sending each $w \in W$ to the unique cluster $x \subseteq W$ containing $w$. Moreover, given a cluster-preserving function $f \colon W \longrightarrow V$, i.e.\ such that
$$w \sim w' ~~~~~\implies~~~~~ f(w) \sim f(w'),$$
it is induced a function $\Refl f \colon \Refl W \longrightarrow \Refl V$, sending $x \in \Refl W$ to the unique cluster $y \in \Refl V$ such that $f(x) \subseteq y$. Let $h \colon X \longrightarrow Y$ be a function between posets. We say that $x \in X$ is \emph{$h$-maximal} if for any $x' \in X$
$$x \prec x' \text{ and } h(x) = h(x') ~~~~~\implies~~~~~ x = x'.$$
\proposition\label{prop:refl}
Let $W$ and $V$ be locally finite preorders. A cluster-preserving function $f \colon W \longrightarrow V$ is a p-morphism if and only if the following conditions hold:
\begin{enumerate}
    \item $\Refl f \colon \Refl W \longrightarrow \Refl V$ is a p-morphism;
    \item if $x \in \Refl W$ is $(\Refl f)$-maximal, %in $(\Refl f)^{-1}((\Refl f)(x))$. (i.e.\ for any $x' \in \Refl W$
    %$$x \prec x' \text{ and } (\Refl f)(x) = (\Refl f)(x') \implies x = x'\text{)}$$
    then the function $x \longrightarrow (\Refl f)(x)$ (restriction of $f$) is surjective.
\end{enumerate}
If $f$ is a p-morphism, then $f$ is surjective if and only if $\Refl f$ is surjective.
\endproposition
\proof
Assume that $f$ is a p-morphism. It is straightforward to check that $\Refl f$ is a p-morphism. We prove condition 2. Let $x \in \Refl W$ be $(\Refl f)$-maximal; %in $(\Refl f)^{-1}((\Refl f)(x))$;
we want to prove that $x \longrightarrow (\Refl f)(x)$ is surjective. For, pick $v' \in (\Refl f)(x)$; if we fix any $w \in x$, we have that $f(w) \prec v'$ (this is because $f(w), v' \in (\Refl f)(x)$ and $(\Refl f)(x)$ is a cluster). Using that $f$ is a p-morphism, we can find $w' \in W$ such that $w \prec w'$ and $f(w') = v'$. Denoting by $x'$ the cluster in $W$ containing $w'$, we have that $x \prec x'$ and $(\Refl f)(x') = (\Refl f)(x)$. By $(\Refl f)$-maximality of $x$, it must be $x = x'$; we then have $w' \in x$ and $f(w') = v'$, witnessing surjectivity of $x \longrightarrow (\Refl f)(x)$.

Vice versa, assume conditions 1. and 2. We have to prove that $f$ is a p-morphism. We start by stability. Let $w, w' \in W$ be such that $w \prec w'$; denote by $x$ and $x'$ the clusters in $W$ containing $w$ and $w'$, respectively. We have that $x \prec x'$; using that $\Refl f$ is a p-morphism, it must be $(\Refl f)(x) \prec (\Refl f)(x')$, implying that $f(w) \prec f(w')$ (since $f(w) \in (\Refl f)(x)$ and $f(w') \in (\Refl f)(x')$). To prove openness, let $w \in W$ and $v' \in V$ be such that $f(w) \prec v'$. Denote by $x$ the cluster in $W$ containing $w$ and by $y'$ the cluster in $V$ containing $v'$; we then have that $(\Refl f)(x) \prec y'$. Using that $\Refl f$ is a p-morphism, we can find $x' \in \Refl W$ such that $x \prec x'$ and $(\Refl f)(x') = y'$. We can assume that $x'$ is $(\Refl f)$-maximal, % in $(\Refl f)^{-1}((\Refl f)(x'))$,
because $\Refl W$ is locally finite, $W$ being so. By condition 2. $x' \longrightarrow (\Refl f)(x')$ is surjective. Since $v' \in y' = (\Refl f)(x')$, we can find $w' \in x'$, hence $w \prec w'$ (by $x \prec x'$), such that $f(w') = v'$.

For the second claim, we need to prove that $f$ is surjective if $\Refl f$ is so (the other direction being straightforward), assuming that $f$ is a p-morphism. Given $v \in V$, let $y$ be the cluster in $V$ containing $v$. By surjectivity of $\Refl f$, we can find $x \in \Refl W$ such that $(\Refl f)(x) = y$. As before, we can assume that $x$ is $(\Refl f)$-maximal. % in $(\Refl f)^{-1}((\Refl f)(x))$, because $\Refl W$ is locally finite.
By what we proved above, $x \longrightarrow (\Refl f)(x) = y$ is surjective, hence we can find $w \in x \subseteq W$ such that $f(w) = v$.
\endproof

The following techniques are taken from \cite[Chapter 6]{GZ}. Given $W_1, \dots, W_n, W$ finite preorders, we say that
$$\textit{\emph{‘‘}$W_1$ and $\cdots$ and $W_n$ yield $W$\emph{"}}$$
if, for each $L \in \NExt \textsf{S4}$ such that $(L\KFr_\fin)^\op$ satisfies the (AP),
$$W_1, \dots, W_n \in L\KFr_\fin ~~~~~\implies~~~~~ W \in L\KFr_\fin.$$

We mention the following notable examples of finite preorders, to which we will apply the notion just introduced.
\begin{enumerate}
    \item[$\bullet$] $C_n = \{1, \dots, n\}$ is the one-cluster frame with $n$ elements ($i \prec j$ for all $i, j \in C_n$);
    \item[$\bullet$] $[n] = \{1, \dots, n\}$ is the $n$-elements linearly ordered chain, with the reverse order ($i \prec j$ if and only if $j \leq i$, for all $i, j \in [n]$, so that $\depth(i) = i$);
    \item[$\bullet$] $D$ is the four elements \say{diamond} poset
    \item[$\bullet$] $W^+$ is the Kripke frame obtained from $W$ by adding a reflexive root to it (an element $\rho$ such that $\rho \prec \rho$ and $\rho \prec w$ for all $w \in W$);
    \item[$\bullet$] $W^-$ is the Kripke frame obtained from $W$ by adding a reflexive final element to it (an element $\varepsilon$ such that $\varepsilon \prec \varepsilon$ and $w \prec \varepsilon$ for all $w \in W$).
    \item[$\bullet$] $\sum_i W_i$ is the ordinary disjoint union;%, i.e.\ the Kripke frame whose underlying set is $\{(w,i) \in \bigcup_{i \in I} W_i \times I\ \vert\ w \in W_i\}$ and such that $(w,i) \prec (w',i')$ if and only if $i = i'$ and $w \prec w'$ in $W_i$;
    \item[$\bullet$] $k W$ is the disjoint union of $k$ copies of $W$ (it coincides with the empty frame for $k = 0$).
\end{enumerate}
For example, the frames $[2]$, $(2[1])^+$, $C_2^-$ and $D$ appear as follows
\[
    \begin{tikzpicture}
    \filldraw (0,0.2) circle (0.1cm);
    \filldraw (0,-0.2) circle (0.1cm);
    \draw (0,-0.1) -- (0,0.1);
    \begin{scope}[shift={(2.5,0)}]
    \filldraw (-0.2,0.2) circle (0.1cm);
    \filldraw (0.2,0.2) circle (0.1cm);
    \filldraw (0,-0.2) circle (0.1cm);
    \draw (0,-0.2) -- (-0.2,0.2) (0,-0.2) -- (0.2,0.2);
    \begin{scope}[shift={(2.5,0)}]
    \filldraw (0,0.5) circle (0.1cm);
    \draw (0,-0.2) circle (0.4cm);
    \filldraw (-0.15,-0.2) circle (0.1cm);
    \filldraw (0.15,-0.2) circle (0.1cm);
    \draw (0,0.4) -- (0,0.2);
    \begin{scope}[shift={(2.5,0)}]
    \filldraw (0,0.4) circle (0.1cm);
    \filldraw (-0.2,0) circle (0.1cm);
    \filldraw (0.2,0) circle (0.1cm);
    \filldraw (0,-0.4) circle (0.1cm);
    \draw (0,-0.4) -- (-0.2,0) (0,-0.4) -- (0.2,0) (0,0.4) -- (-0.2,0) (0,0.4) -- (0.2,0);
    %\draw (0,0.4) circle (0.4cm);
    %\filldraw (-0.15,0.4) circle (0.1cm);
    %\filldraw (0.15,0.4) circle (0.1cm);
    %\filldraw (-0.2,-0.3) circle (0.1cm);
    %\filldraw (0.2,-0.3) circle (0.1cm);
    %\filldraw (0,-0.7) circle (0.1cm);
    %\draw (-0.05,-0.65) -- (-0.2,-0.4) (0.05,-0.65) -- (0.2,-0.4) (-0.2,-0.2) -- (-0.1,0.01) (0.2,-0.2) -- (0.1,0.01);
    \end{scope}
    \end{scope}
    \end{scope}
    \end{tikzpicture}
\]

\lemma\label{lem:yield}
Let $W$ be a finite preorder.
\begin{enumerate}
    \item[(o)] if $W$ has an external (resp. internal) cluster of cardinality $\geq n$, then $W$ yield $C_n$ (resp. $C_n^-$);\footnote{Recall that a cluster $C \subseteq W$ is said to be external if $\depth(w) = 1$ for some $w \in C$ (i.e.\ for all $w \in C$), where $\depth(w)$ is the maximum cardinality of sequences $w = w_1 \prec w_2 \prec \cdots \prec w_n$, with $w_i \neq w_{i+1}$ and $w_{i+1} \nprec w_i$. Otherwise, the cluster $C$ is said to be internal.}
    \item[(i)] $W$ and $C_n$ yield $W^e_n$, where $W^e_n$ is obtained from $W$ by replacing all external clusters in it by $C_n$;
    \item[(ii)] $W$ and $C_n^-$ yield $W^p_n$, where $\{p\} \subseteq W$ is a singleton internal cluster and $W^p_n$ is obtained from $W$ by replacing the cluster $\{p\}$ by $C_n$;
    \item[(iii)] $C_3$ and $C_n$ yield $C_{n+1}$;
    \item[(iv)] $C_3^-$ and $C_n^-$ yield $C_{n+1}^-$;
    \item[(v)] $C_2^{++}$ and $C_2^{++}$ yield $D^e_2$;
    \item[(vi)] $(C_2^-)^{++}$ and $(C_2^-)^{++}$ yield $(D^e_2)^-$.
\end{enumerate}
\endlemma
\proof
Throughout this proof, fix $L \in \NExt \textsf{S4}$ such that $(L\KFr_\fin)^\op$ satisfies the (AP).

(o) is an easy consequence of the fact that $L\KFr_\fin$ is closed under p-morphic images (this is because taking p-morphic images preserves the logic, see \cite{misha}).\footnote{Here, the assumption that $(L\KFr_\fin)^\op$ satisfies the (AP) does not play any role.}

To prove (i) and (ii) we use Proposition~\ref{prop:refl}. We show (ii) ((i) can be proved in a similar way). Assume that $W, C_n^- \in L\KFr_\fin$ and let $\{p\} \subseteq W$ be a singleton internal cluster. Consider the following surjective p-morphisms: the first one is $f \colon W \longrightarrow [2]$, sending all the points $w \in W$ such that $w \succ p$ and $w \neq p$ to $1$, and all the other points to $2$; the second one is $g \colon C_n^- \longrightarrow [2]$ sending the final point to $1$, and all the points in $C_n$ to $2$. Since $L\KFr_\fin$ has the dual (AP), we can find surjective p-morphisms $h \colon A \longrightarrow W$ and $k \colon A \longrightarrow C_n^-$ in $L\KFr_\fin$ such that $fh = gk$. To prove that $W_n^p \in L\KFr_\fin$, we define a surjective p-morphism $\hat{h} \colon A \longrightarrow W_n^p$ (we use again logic preservation of p-morphic images). For, observe that $\Refl W_n^p = \Refl W$ and that we have a surjective p-morphism $\pi_n^p \colon W_n^p \longrightarrow W$ sending all the points in the new cluster to $p$. By Proposition~\ref{prop:refl}, the p-morphism $h \colon A \longrightarrow W$ can be decomposed as $\hat{h} \colon A \longrightarrow W_n^p$ followed by $\pi_n^p \colon W_n^p \longrightarrow W$, provided that:
\begin{gather*}
\text{($*$) for each $(\Refl h)$-maximal cluster $C \in \Refl A$ such that $(\Refl h)(C) = \{p\}$,}\\
\text{there exists a surjective function $C \longrightarrow C_n$.}
\end{gather*}
To see that ($*$) holds, observe that, if $C$ is $(\Refl h)$-maximal and $(\Refl h)(C) = \{p\}$, then $C$ must be $(\Refl k)$-maximal: this is because $(\Refl k)(C) = C_n$ (since $(\Refl gk)(C) = (\Refl fh)(C) = \{2\}$) and for any other $C' \in \Refl A$, if $C \prec C'$ and $(\Refl k)(C') = C_n$, then $(\Refl h)(C') = \{p\}$ (since $\{p\} = (\Refl h)(C) \prec (\Refl h)(C')$ and $(\Refl fh)(C') = (\Refl gk)(C') = \{2\}$). By Proposition~\ref{prop:refl}, the restriction $C \longrightarrow (\Refl k)(C) = C_n$ of $k$ must be surjective. As we said, this defines a p-morphism $\hat{h} \colon A \longrightarrow W_n^p$; the surjectivity of $\hat{h}$ now follows from that of $h$, since $\Refl \hat{h} = \Refl h$ (again, by Proposition~\ref{prop:refl}).
%consider the (unique) surjective p-morphisms $W \longrightarrow [1]$ and $C_n \longrightarrow [1]$ and a pair of surjective p-morphisms $f \colon U \longrightarrow W$ and $g \colon U \longrightarrow C_n$ (which necessarily coamalgamates the starting pair, $[1]$ being the terminal object in $S4\KFr_\lf$). It is straightforward to see that there is a surjective p-morphism $W^e_n \longrightarrow W$ and that $f$ factorizes trough it ($g$ forces the external clusters of $U$ to be of cardinality greater or equal than $n$); moreover, $U \longrightarrow W^e_n$ of such factorization is surjective.

To prove (iii)-(vi) we make use of the following notion: say that a pair of surjective p-morphisms $f \colon W \longrightarrow V$ and $g \colon U \longrightarrow V$ is \emph{thin} %\footnote{See \cite[Exercise 6.4.1]{GZ}}
if, for each $v \in V$,
$$|f^{-1}(v)| > 1 ~~~~~\implies~~~~~ |g^{-1}(v)| = 1$$
(and the same exchanging the roles of $f$ and $g$). If the pair $f$ and $g$ as above is thin, then $W$ and $U$ yield $P_{f,g}$ (recall that $P_{f,g}$ is the $\DGrph$-pullback of $f$ and $g$, which is given by the set-theoretic pullback equipped with the restricted product relation). This is because, for each pair $h \colon A \longrightarrow W$ and $k \colon A \longrightarrow U$ of surjective p-morphisms such that $fh = gk$, the unique stable function $A \longrightarrow P_{f,g}$ sending $a \in A$ to $(h(a),k(a)) \in P_{f,g}$ is a surjective p-morphism; this ensures that $P_{f,g} \in L\KFr_\fin$ whenever $L\KFr_\fin$ contains a coamalgam of $f$ and $g$, since taking p-morphic images preserves the logic. With this fact at hand, we show (v) (the others can be proved in a similar way). For, the pair $f \colon C_2^{++} \longrightarrow C_2^+$ and $g \colon C_2^{++} \longrightarrow C_2^+$ depicted below
\[
    \begin{tikzpicture}
    \draw (0,0.4) circle (0.4cm);
    \filldraw[red] (-0.15,0.4) circle (0.1cm);
    \filldraw[blue] (0.15,0.4) circle (0.1cm);
    \filldraw[red] (0,-0.3) circle (0.1cm);
    \filldraw (0,-0.7) circle (0.1cm);
    \draw (0,-0.6) -- (0,-0.4) (0,-0.2) -- (0,0);
    \begin{scope}[shift={(2,0)}]
    \draw[->] (-0.8,0) -- (0.8,0) node [midway, above=1pt, fill=white] {$f$};
    \begin{scope}[shift={(2,0)}]
    \draw (0,0.2) circle (0.4cm);
    \filldraw[red] (-0.15,0.2) circle (0.1cm);
    \filldraw[blue] (0.15,0.2) circle (0.1cm);
    \filldraw (0,-0.5) circle (0.1cm);
    \draw (0,-0.4) -- (0,-0.2);
    \begin{scope}[shift={(2,0)}]
    \draw[<-] (-0.8,0) -- (0.8,0) node [midway, above=1pt, fill=white] {$g$};
    \begin{scope}[shift={(2,0)}]
    \draw (0,0.4) circle (0.4cm);
    \filldraw[red] (-0.15,0.4) circle (0.1cm);
    \filldraw[blue] (0.15,0.4) circle (0.1cm);
    \filldraw[blue] (0,-0.3) circle (0.1cm);
    \filldraw (0,-0.7) circle (0.1cm);
    \draw (0,-0.6) -- (0,-0.4) (0,-0.2) -- (0,0);
    \end{scope}
    \end{scope}
    \end{scope}
    \end{scope}
    \end{tikzpicture}
\]
is obviously thin, and $P_{f,g}$ is
\[
    \begin{tikzpicture}
    \draw (0,0.4) circle (0.4cm);
    \filldraw (-0.15,0.4) circle (0.1cm);
    \filldraw (0.15,0.4) circle (0.1cm);
    \filldraw (-0.2,-0.3) circle (0.1cm);
    \filldraw (0.2,-0.3) circle (0.1cm);
    \filldraw (0,-0.7) circle (0.1cm);
    \draw (0,-0.7) -- (-0.2,-0.3) (0,-0.7) -- (0.2,-0.3) (-0.2,-0.3) -- (-0.1,0.01) (0.2,-0.3) -- (0.1,0.01);
    \end{tikzpicture}
\]
Then $C_2^{++}$ and $C_2^{++}$ yield $D_2^e$.
\endproof

To state the key result of this section, we complete the list of invariants attached to a locally finite preorder $W$ introduced above Proposition~\ref{prop:regnecpos}. Let $\exter(W)$ be the maximum cardinality of external clusters in $W$ and let $\inter(W)$ be the maximum cardinality of internal clusters in $W$. As before, $\exter(W) = \omega$ if the maximum does not exist and $\exter(W) = 0$ if $W=\emptyset$ (and the same for $\inter$).
\proposition\label{prop:regnec}
Let $L$ be a normal modal logic with the finite model property and assume that $L \in \NExt \textup{\textsf{S4}}$. If $L\MA_\fin$ satisfies the \emph{(AP)}, then $L = \textup{\textsf{Grz.3}}$, or $L\KFr_\fin$ is one of the following:
\begin{enumerate}
    \item[(I)] $\{W \in \textup{\textsf{S4}}\KFr_\fin\ \vert\ \depth(W) \leq l,\ \width(W) \leq k,\ \exter(W) \leq m,\ \inter(W) \leq n\}$, for some $l \in \{0,1,2\}$ and some $k,m,n \in \{1,2,\omega\}$;
    \item[(II)] $\{W \in \textup{\textsf{S4}}\KFr_\fin\ \vert\ \numberext(W) \leq s,\ \exter(W) \leq m,\ \inter(W) \leq n\}$, for some $s \in \{1,\omega\}$ and some $m,n \in \{1,2,\omega\}$.
\end{enumerate}
\begin{comment}
\begin{enumerate}
    \item $L' + \textit{be}_m + \textit{bi}_n$, where $m, n \in \{1,2,\omega\}$ and $L'$ is one between
    \begin{enumerate}
        \item \textup{\textsf{S4}};
        \item $\textup{\textsf{S4.2}}$;% = \textup{\textsf{S4}} + \textit{bf}_1$;
        \item $\textup{\textsf{S4}} + \textit{bd}_2$;
        \item $\textup{\textsf{S4}} + \textit{bd}_2 + \textit{bw}_2$;
        \item $\textup{\textsf{S4.2}} + \textit{bd}_2$;
    \end{enumerate}
    \item \textup{\textsf{Grz.3}};
    \item $\textup{\textsf{S5}} + \textit{be}_m$, where $m \in \{1,2,\omega\}$.
\end{enumerate}
\end{comment}
\endproposition
\proof
%Let $L \in \NExt \textsf{S4}$ and let $\rho L$ be the so-called \emph{superintuitionistic fragment} of $L$ (see \cite[Chapter 9.6]{MR3363829}). We denote by $\rho L\HA$ the variety of Heyting algebras validating $\rho L$, and by $\rho L\HA_\fin$ the class of the finite members of $\rho L \HA$. Assume that $L\MA_\fin$ has the (AP). By \cite[Proposition 8.3]{MR2153890} restricted to finite algebras (see \cite[\S 3, Theorem 3.1]{MR163850}) $\rho L\HA_\fin$ has the (AP). %By \cite[Theorem 6.39]{MR2153890}, also $\rho L\HA$ has the (AP).
%By Maksimova's characterization of amalgamable varieties of Heyting algebras \cite[Theorem 2]{Maksimova1977CraigsTI}, we have that $L' \subseteq L \subseteq L' + \textsf{Grz}$, where $L'$ is one between \textsf{S4}, \textsf{S4.2}, $\textup{\textsf{S4}} + \textit{bd}_2$, $\textup{\textsf{S4}} + \textit{bd}_2 + \textit{bw}_2$, $\textup{\textsf{S4.2}} + \textit{bd}_2$, \textsf{S4.3} and \textsf{S5} (see \cite[Theorem 7]{Maksimova1979InterpolationTI}).
Assume that $L\MA_\fin \simeq (L\KFr_\fin)^\op$ satisfies the (AP). We apply the results of Lemma~\ref{lem:yield}. We have two possible situations: either all the frames in $L\KFr_\fin$ are partial orders, or $L\KFr_\fin$ contains a preorder that is not a partial order. In the latter, either $C_2 \in L\KFr_\fin$ or ($C_2^-$ is there and $C_2$ is not), by (o).
\begin{enumerate}
    \item In the former case, $L\KFr_\fin$
    %\supseteq \{W \in \textsf{S4}\KFr_\fin\ \vert\ \Refl W \in L\KFr_\fin,\ \delta^e(W) \leq 2,\ \delta^i(W) \leq 1\}$,
    contains all the finite preorders $W$ such that $\Refl W \in L\KFr_\fin$, $\delta^e(W) \leq 2$ and $\delta^i(W) \leq 1$,
    by (i). If there is more, then $C_3 \in L\KFr_\fin$ or ($C_2^-$ is there and $C_3$ is not), by (o).
    \begin{enumerate}
        \item In the former case, $C_n \in L\KFr_\fin$ for all $n \geq 1$ (by (iii)); as a consequence, $L\KFr_\fin$ contains all $W$ such that $\Refl W \in L\KFr_\fin$ and $\delta^i(W) \leq 1$, by (i). If we have more, then $C_2^- \in L\KFr_\fin$ (by (o)), hence $L\KFr_\fin$ contains all $W$ such that $\Refl W \in L\KFr_\fin$ and $ \delta^i(W) \leq 2$, by (ii). If we have even more, then $C_3^- \in L\KFr_\fin$ (by (o)), which implies $C_n^- \in L\KFr_\fin$ for all $n \geq 1$ (by (iv)); as a consequence, $L\KFr_\fin$ contains all $W$ such that $\Refl W \in L\KFr_\fin$, by (ii).
        \item In the latter case, $L\KFr_\fin$ contains all $W$ such that $\Refl W \in L\KFr_\fin$, $\delta^e(W) \leq 2$ and $\delta^i(W) \leq 2$, by (ii). If we have more, then $C_3^- \in L\KFr_\fin$ (by (o)). Reasoning as above, $C_n^- \in L\KFr_\fin$ for all $n \geq 1$; as a consequence, $L\KFr_\fin$ contains all $W$ such that $\Refl W \in L\KFr_\fin$ and $\delta^e(W) \leq 2$, by (ii).
    \end{enumerate}
    \item In the latter case, $L\KFr_\fin$ contains all $W$ such that $\Refl W \in L\KFr_\fin$, $\delta^e(W) \leq 1$ and $\delta^i(W) \leq 2$, by (ii). If there is more, then $L\KFr_\fin$ contains all $W$ such that $\Refl W \in L\KFr_\fin$ and $\delta^e(W) \leq 1$, reasoning as before.
\end{enumerate}
Consider now $L_0 \coloneqq L + \textsf{Grz} \in \NExt \textsf{Grz}$; observe that $L_0\KFr_\fin = \{\Refl W\ \vert\ W \in L\KFr_\fin\}$ (this is because the projection $W \longrightarrow \Refl W$ is a surjective p-morphism, and taking p-morphic images is logic preserving \cite{misha}). It is straightforward to prove that $(L_0\KFr_\fin)^\op$ satisfies the (AP), because $(L\KFr_\fin)^\op$ does so, by hypothesis. This means that $L_0$ is one of the $8$ logics of Proposition~\ref{prop:regnecpos}. Putting everything together, we have that
$$L\KFr_\fin = \{W \in \textsf{S4}\KFr_\fin\ \vert\ \Refl W \in L_0\KFr_\fin,\ \exter(W) \leq m,\ \inter(W) \leq n\},$$
for some $m,n \in \{1,2,\omega\}$, where $L_0$ is one of the logics of Proposition~\ref{prop:regnecpos}. To conclude, we show that, in case $L_0 = \textsf{Grz.3}$, then the argument above trivializes, forcing $L = \textsf{Grz.3}$.
\begin{enumerate}
    \item If $C_2 \in L\KFr_\fin$, then $C_2^{++} \in L\KFr_\fin$ (by (i)), since $[3] \in \textsf{Grz.3}\KFr_\fin \subseteq L\KFr_\fin$. As a consequence, $D^e_2 \in L\KFr_\fin$ (by (v)), hence $D = \Refl D^e_2 \in \textsf{Grz.3}\KFr_\fin$.
    \item If $C_2^- \in L\KFr_\fin$, then $(C_2^-)^{++} \in L\KFr_\fin$ (by (ii)), since $[4] \in \textsf{Grz.3}\KFr_\fin \subseteq L\KFr_\fin$. As a consequence, $(D^e_2)^- \in L\KFr_\fin$ (by (vi)), hence $D^- = \Refl (D^e_2)^- \in \textsf{Grz.3}\KFr_\fin$.
\end{enumerate}
Both cases lead to a contradiction, proving that all the preorders in $L\KFr_\fin$ are partial orders.
\endproof

\remark
Observe that the above classification subsumes Maksimova's (\cite[Theorem 7]{Maksimova1979InterpolationTI}): by Proposition~\ref{prop:APalg}, if $L \in \NExt \textup{\textsf{K4}}$ has the finite model property, then amalgamability of $L\MA$ implies amalgamability of $L\MA_\fin$. However, as we will see at the end of the section, the converse does not hold.
\endremark

\remark
Some of the logics mentioned in the Table in Section~\ref{sec:profinite} appear in Proposition~\ref{prop:regnec}. For example, \textsf{S5} is described in (I), while \textsf{Grz}, \textsf{S4} and \textsf{S4.2} are described in (II).
\endremark

\subsection{Completeness of the classification}
We now want to prove that the list of Proposition~\ref{prop:regnec} provides a complete classification of the amalgamability of $L\MA_\fin \simeq (L\KFr_\fin)^\op$, which is equivalent to the regularity of $L\KFr_\lf \simeq (\Pro L\MA_\fin)^\op$, by Theorem~\ref{teo:regular} and Proposition~\ref{prop:APfin}.

We first consider the logics corresponding to case (II) of Proposition~\ref{prop:regnec}, namely
$$L\KFr_\fin = \{W \in \textup{\textsf{S4}}\KFr_\fin\ \vert\ \numberext(W) \leq s,\ \exter(W) \leq m,\ \inter(W) \leq n\},$$
for some $s \in \{1,\omega\}$ and some $m,n \in \{1,2,\omega\}$. Observe that $\numberext(W) \leq 1$ if and only if $W$ is locally confluent, meaning that $w^*$ has a maximum for each $w \in W$.
\proposition\label{prop:preoconfl}
If $L$ is as in (II) of \textup{Proposition~\ref{prop:regnec}}, then $L\MA_\fin$ satisfies the \textup{(AP)}.
\endproposition
\proof
Let $L$ be as in the hypothesis and let $f \colon W \longrightarrow V$ and $g \colon U \longrightarrow V$ be surjective p-morphisms in $L\KFr_\fin$. Consider the subset of the set-theoretic product $\Refl W \times \Refl U$
$$A \coloneqq \{(x,y)\ \vert\ (\Refl f)(x) = (\Refl g)(y) \text{ and } \exists w \in x, u \in y \text{ s.t.\ } f(w) = g(u)\},$$
and say that
$$(x,y) \prec (x',y') ~~~~~\iff~~~~~ x \prec x' \text{ and } y \prec y';$$
it is straightforward to verify that $A \in (L+\textsf{Grz})\KFr_\fin$: it is a finite poset and $\numberext(A) \leq 1$ if both $\numberext(W) \leq 1$ and $\numberext(U) \leq 1$. We then have a diagram
\[\begin{tikzcd}
	A & {\Refl U} \\
	{\Refl W} & {\Refl V}
	\arrow["{\pi_U}", from=1-1, to=1-2]
	\arrow["{\pi_W}"', from=1-1, to=2-1]
	\arrow["{{\Refl g}}", from=1-2, to=2-2]
	\arrow["{{\Refl f}}"', from=2-1, to=2-2]
\end{tikzcd}\]
in $(L+\textsf{Grz})\KFr_\fin$. The projections are easily checked to be surjective, hence such a diagram is a coamalgam in $(L+\textsf{Grz})\KFr_\fin$. We now \textit{blow-up} $A$ to a coamalgam in $L\KFr_\fin$ for the original diagram. For each $(x,y) \in A$, we define a set $C(x,y)$ and functions $C(x,y) \longrightarrow x$ and $C(x,y) \longrightarrow y$ making the diagram
\[\begin{tikzcd}
	{C(x,y)} & y \\
	x & {(\Refl f)(x) = (\Refl g)(y)}
	\arrow[from=1-1, to=1-2]
	\arrow[from=1-1, to=2-1]
	\arrow[from=1-2, to=2-2]
	\arrow[from=2-1, to=2-2]
\end{tikzcd}\]
commute. We distinguish the following cases.
\begin{enumerate}
    \item[(i)]  Assume first that $x$ is $(\Refl f)$-maximal %in $(\Refl f)^{-1}((\Refl f)(x))$
    and that $y$ is $(\Refl g)$-maximal. %in $(\Refl g)^{-1}((\Refl g)(y))$.
    By Proposition~\ref{prop:refl}, this means that $x \longrightarrow (\Refl f)(x)$ and $y \longrightarrow (\Refl g)(y)$ are both surjective. Consider $v \in (\Refl f)(x) = (\Refl g)(y)$; the inverse images $x_v \coloneqq \{w \in x\ \vert\ f(w)=v\}$ and $y_v \coloneqq \{u \in y\ \vert\ g(u)=v\}$ are both non-empty. We can assume, without loss of generality, that the cardinality of $x_v$ is bigger or equal to the cardinality of $y_v$ (if this is not the case, exchange the roles of $x$ and $y$ in what follows). Let $z_v = x_v$, with inclusion $z_v \longrightarrow x$, and let $z_v \longrightarrow y$ be any function mapping surjectively to $y_v$. If we let $C(x,y)$ be the disjoint union of the $z_v$'s, varying $v \in (\Refl f)(x) = (\Refl g)(y)$, we obtain surjective functions $C(x,y) \longrightarrow x$ and $C(x,y) \longrightarrow y$ making the diagram above commute. Observe that, if $|x| \leq 2$ and $|y| \leq 2$, then $|C(x,y)| \leq 2$ (and the same if the bound is $1$).% This case proves the dual (AP) for $(\mathsf{S5} + \textit{be}_k)\KFr_\fin$, with $k \in \{1,2,\omega\}$ (check that $|C(x,y)| \leq k$, if $|x| \leq k$ and $|y| \leq k$).
    
    \item[(ii)] Assume that $y$ is not $(\Refl g)$-maximal. %in $(\Refl g)^{-1}((\Refl g)(y))$.
    Let
    $$C(x,y) \coloneqq \{u \in y\ \vert\ \exists w \in x \text{ s.t. } f(w) = g(u)\},$$
    with inclusion $C(x,y) \longrightarrow y$, and let $C(x,y) \longrightarrow x$ choose for each $u \in C(x,y)$ some $w \in x$ such that $f(w) = g(u)$. Observe that, in case $x$ is $(\Refl f)$-maximal, %in $(\Refl f)^{-1}((\Refl f)(x))$,
    then $x \longrightarrow (\Refl f)(x)$ is surjective, hence $C(x,y) = y$ and $C(x,y) \longrightarrow y$ is the identity.
    
    \item[(iii)] Finally, assume that $y$ is $(\Refl g)$-maximal %in $(\Refl g)^{-1}((\Refl g)(y))$
    and that $x$ is not $(\Refl f)$-maximal. %in $(\Refl f)^{-1}((\Refl f)(x))$.
    In this case, the construction is analogous to the previous one, exchanging the roles of $x$ and $y$: let $C(x,y) = x$, with identity $C(x,y) \longrightarrow x$, and let $C(x,y) \longrightarrow y$ choose for each $w \in C(x,y) = x$ some $u \in y$ such that $f(w) = g(u)$ (it is well defined because $y \longrightarrow (\Refl g)(y)$ is surjective, by Proposition~\ref{prop:refl}).
\end{enumerate}
Call $\hat{A}$ the preorder obtained by replacing the point $(x,y)$ in $A$ with the cluster $C(x,y)$, defined as above: in other words,
$$\hat{A} = \{(x,y,a)\ \vert\ (x,y) \in A,\ a \in C(x,y)\}$$
and, for each $(x,y,a), (x',y',a') \in \hat{A}$,
$$(x,y,a) \prec (x',y',a') ~~~~~\iff~~~~~ x \prec x' \text{ and } y \prec y'.$$
We have that $\Refl \hat{A} = A$ (we are identifying $C(x,y)$ and $(x,y)$). Additionally, the construction of the $C(x,y)$'s gives two cluster-preserving functions $\Pi_W \colon \hat{A} \longrightarrow W$ and $\Pi_U \colon \hat{A} \longrightarrow U$, such that $\Refl \Pi_W = \pi_W$ and $\Refl \Pi_U = \pi_U$, satisfying $f \Pi_W = g \Pi_U$. To prove that the pair $\Pi_W$ and $\Pi_U$ is a coamalgam for $f$ and $g$ in $L\KFr_\fin$, we need to prove that: (a) $\hat{A}$ belongs to $L\KFr_\fin$; (b) $\Pi_W$ and $\Pi_U$ are p-morphisms (by Proposition~\ref{prop:refl}, if that is the case, then $\Pi_W$ and $\Pi_U$ are surjective, because $\Refl \Pi_W = \pi_W$ and $\Refl \Pi_U = \pi_U$ are so).

To prove (a) we only need to check that $\exter(\hat{A}) \leq m$ and that $\inter(\hat{A}) \leq n$. Observe that, for each $k \in \{1,2,\omega\}$, if $|x| \leq k$ and $|y| \leq k$, then $|C(x,y)| \leq k$: we already observed this for case (i); for cases (ii) and (iii), it follows from the fact that $C(x,y)$ is a subset of either $x$ or $y$. First, assume that $C(x,y)$ is external in $\hat{A}$; this is equivalent to say that $(x,y)$ is external in $A$, i.e.\ that $x$ is external in $W$ and $y$ is external in $U$. In this situation, $|x| \leq m$ and $|y| \leq m$, hence $|C(c,y)| \leq m$. Finally, if $C(x,y)$ is internal in $\hat{A}$, then at least one between $x$ and $y$ is internal: if both of them are internal, then $|x| \leq n$ and $|y| \leq n$, hence $|C(x,y)| \leq n$; if $x$ is external and $y$ is internal, then we are in case (ii), hence $|C(x,y)| \leq n$, because $C(x,y) \subseteq y$ and $|y| \leq n$.

To prove (b) we use Proposition~\ref{prop:refl}. We already know that $\Refl \Pi_W = \pi_W$ is a p-morphism. To show the other condition, we first prove that, given $(x,y) \in \Refl \hat{A} = A$,
\begin{align*}
(x,y) \text{ is } (\Refl \Pi_W)\text{-maximal} ~~~~~\implies~~~~~ y \text{ is } (\Refl g)\text{-maximal}.\tag{$*$}
\end{align*}
For, let $(x,y) \in A$ be $(\Refl \Pi_W)$-maximal. Pick $y' \in \Refl U$ $(\Refl g)$-maximal such that $y \prec y'$ and $(\Refl g)(y') = (\Refl g)(y)$ (it exists because $\Refl U$ is finite). Then $(x,y') \in A$ (because $y' \longrightarrow (\Refl g)(y') = (\Refl f)(x)$ is surjective), $(x,y) \prec (x,y')$ and $(\Refl \Pi_W)(x,y) = x = (\Refl \Pi_W)(x,y')$. By maximality, we have that $(x,y) = (x,y')$; in particular, $y=y'$ is $(\Refl g)$-maximal. This proves ($*$). We are now ready to prove condition 2. of Proposition~\ref{prop:refl} for $\Pi_W$. Let $(x,y) \in \Refl \hat{A} = A$ be $(\Refl \Pi_W)$-maximal. By ($*$), $y$ is $(\Refl g)$-maximal. We are then either in case (i), or in case (iii); in both cases, $C(x,y) \longrightarrow x = (\Refl \Pi_W)(C(x,y))$ is surjective. This proves that $\Pi_W$ is a p-morphism; similarly, it is possible to prove that $\Pi_U$ is a p-morphism.
\endproof

We now consider the case $L = \textsf{Grz.3}$. Observe that $W \in \textsf{Grz.3}\KFr_\fin$ if and only if $W$ is a finite locally linear poset, where being locally linear means that $w^*$ is a linear order for each $w \in W$.

Consider $f \colon W \longrightarrow V$ and $g \colon U \longrightarrow V$ in $\textsf{Grz.3}\KFr_\fin$. Recall that $P_{f,g}$ denotes the $\DGrph$-pullback of $f$ and $g$. The pullback in $\textsf{Grz.3}\KFr_\fin$ can be described as follows (see \cite{Godel}). Let $\Chains(P_{f,g})$ be the set of \emph{linearly ordered bisimulations} in $P_{f,g}$, namely the linearly ordered subsets $C \subseteq P_{f,g}$ satisfying the following property: for each $(w,u) \in C$ and $w' \in W$, if $w \prec w'$, then there exists $u' \in U$ such that $u \prec u'$ and $(w',u') \in C$ (and symmetrically). Moreover, for $C, D \in \Chains(P_{f,g})$, we say that $C \prec D$ if and only if: (i) $D \subseteq C$; (ii) for each $(w,u) \in D$ and $(w',u') \in C$, if $(w,u) \prec (w',u')$ in $P_{f,g}$, then $(w',u') \in D$. It is straightforward to prove that $\Chains(P_{f,g}) \in \textsf{Grz.3}\KFr_\fin$. We can also define projections $p_W \colon \Chains(P_{f,g}) \longrightarrow W$ and $p_U \colon \Chains(P_{f,g}) \longrightarrow U$: if the minimum of $C \in \Chains(P_{f,g})$ is $(w_0,u_0)$, then $p_W(C) = w_0$ and $p_U(C) = u_0$. The maps $p_W$ and $p_U$ are easily seen to be p-morphisms and, obviously, $f p_W = g p_U$.%; moreover, given any pair of p-morphisms $h \colon A \longrightarrow W$ and $k \colon A \longrightarrow U$ in $\textsf{Grz.3}\KFr_\fin$, the p-morphism $A \longrightarrow \Chains(P_{f,g})$ witnessing the universal property of the pullback sends $a \in A$ to the linearly ordered bisimulation $\{(h(a'),k(a'))\ \vert\ a \prec a'\}$.

\proposition\label{prop:poslin}
$\textup{\textsf{Grz.3}}\MA_\fin$ satisfies the \textup{(AP)}.
\endproposition
\proof
With the description of pullbacks in $\textup{\textsf{Grz.3}}\KFr_\fin$ given above, it is sufficient to prove that $p_W$ is surjective whenever $g$ is surjective. This means that, for each pair $f \colon W \longrightarrow V$ and $g \colon U \longrightarrow V$ of p-morphisms in $\textup{\textsf{Grz.3}}\KFr_\fin$, and for each $w_0 \in W$:
\begin{gather*}
\text{$(*)$ if $g$ is surjective, then there exists a linearly ordered bisimulation}\\
\text{$C$ in $P_{f,g}$ such that $p_W(C) = w_0$.}
\end{gather*}
Let $f$, $g$ and $w_0$ as above; we prove $(*)$ by induction on the depth $\depth(w_0)$. Assume that $g$ is surjective. We can fix $u_0 \in U$ such that $v_0 \coloneqq f(w_0) = g(u_0)$; without loss of generality, we assume that $u_0$ is $g$-maximal.\footnote{Recall that this means that for each $u \in U$, $u_0 \prec u$ and $g(u_0) = g(u)$ imply $u_0 = u$.} If $\depth(w_0) = 1$, then $C = \{(w_0,u_0)\}$ is a linearly ordered bisimulation in $P_{f,g}$ and $p_W(C) = w_0$. Now assume that $\depth(w_0) > 1$ and let $w_1$ be the immediate successor of $w_0$ in $W$. Obviously, $\depth(w_1) < \depth(w_0)$; by induction hypothesis applied to the triple $f_0 \colon w_0^* \longrightarrow v_0^*$, $g_0 \colon u_0^* \longrightarrow v_0^*$ (restrictions of $f$ and $g$) and $w_1 \in w_0^*$, we can find a linearly ordered bisimulation $D$ in $P_{f_0,g_0}$ such that $p_{w_0^*}(D) = w_1$. Without loss of generality, we can assume that $D$ has maximal length. But then $C = \{(w_0,u_0)\} \cup D$ is a linearly ordered bisimulation in $P_{f,g}$ and $p_W(C) = w_0$.
\endproof

The following theorem summarizes all the results obtained in this section, together with the syntactic characterization obtained in \cite{infinitary}. We present the logics in terms of their locally finite frames rather than their finite frames; this is because we want to emphasize the regularity of $L\KFr_\lf \simeq (\Pro L\MA_\fin)^\op$, which will be strengthened to Barr exactness in the following section. Indeed, using the fact that generated subframes, p-morphic images and disjoint unions are logic preserving \cite{misha}, one can show that $L\KFr_\lf = \{W \in \KFr_\lf\ \vert\ (\forall w \in W)(w^* \in L\KFr_\fin)\}$.
\theorem\label{teo:regchar}
Let $L$ be a normal modal logic with the finite model property and assume that $L \in \NExt \textup{\textsf{S4}}$. The following are equivalent:
\begin{enumerate}
    \item $L$ satisfies the infinitary Craig’s Interpolation Theorem for the global consequence relation;
    \item $L\KFr_\lf$ is regular;
    \item $\Pro L\MA_\fin$ satisfies the \textup{(AP)};
    \item $L\MA_\fin$ satisfies the \textup{(AP)};
    \item $L = \textup{\textsf{Grz.3}}$, or $L\KFr_\lf$ is one of the following:
    \begin{enumerate}
        \item[(I)] $\{W \in \textup{\textsf{S4}}\KFr_\lf\ \vert\ \depth(W) \leq l,\ \width(W) \leq k,\ \exter(W) \leq m,\ \inter(W) \leq n\}$, for some $l \in \{0,1,2\}$ and some $k,m,n \in \{1,2,\omega\}$;
        \item[(II)] $\{W \in \textup{\textsf{S4}}\KFr_\lf\ \vert\ \numberext(W) \leq s,\ \exter(W) \leq m,\ \inter(W) \leq n\}$, for some $s \in \{1,\omega\}$ and some $m,n \in \{1,2,\omega\}$.
    \end{enumerate}
\end{enumerate}
\endtheorem
\proof
It remains to prove that $L\MA_\fin$ satisfies the (AP) for all the logics described by (I). This follows from Proposition~\ref{prop:APalg}, since $L\MA$ satisfies the (AP) for all such logics, as shown in \cite{Maksimova1980InterpolationTI}.
\endproof

\remark\label{rem:exc}
As we anticipated, our classification differs from that of Maksimova. Namely, $\textsf{Grz.3}\MA_\fin$ satisfies the (AP) (by Proposition~\ref{prop:poslin}), while $\textsf{Grz.3}\MA$ does not (see \cite{Maksimova1982AbsenceOT}).
\endremark

\section{Coexactness}\label{sec:coexactness}
Also Barr exactness has a syntactic counterpart in the infinitary modal calculi presented in \cite{infinitary}: it corresponds to the fact that every formula $\rho(X_0,X_1)$, where $X_0$ and $X_1$ are isomorphic sets of variables, for which the logic proves reflexivity, symmetry, and transitivity, is of the form $\bigwedge_i (\varphi_i(X_0) \leftrightarrow \varphi_i(X_1))$. As we did for regularity, we want to characterize the normal modal logics $L \in \NExt \textsf{S4}$ with the finite model property for which $L\KFr_\lf$ is Barr exact.

Barr exact categories are defined as those regular categories whose equivalence relations are in one-to-one correspondence with quotients.
\definition\label{defn:eq}
Fix an object $W$ in a (regular) category $\mathbf{C}$. An \emph{equivalence relation} over $W$ is a jointly monic pair of parallel morphisms $p_0, p_1 \colon X \longrightarrow W$ (called \emph{projections}) for which
\begin{enumerate}
    \item[r)] there exists $r \colon W \longrightarrow X$ such that $p_i r = \text{id}_W$ for $i = 0, 1$;
    \item[s)] there exists $s \colon X \longrightarrow X$ such that $p_0 s = p_1$ and $p_1 s = p_0$;
    \item[t)] for all $h_0, h_1 \colon T \longrightarrow X$ such that $p_0 h_1 = p_1 h_0$, there exists a (unique, the pair $p_0$, $p_1$ being jointly monic) $h \colon T \longrightarrow X$ such that $p_i h = p_i h_i$ for $i = 0, 1$.
\end{enumerate}
\enddefinition

Given a morphism $f \colon W \longrightarrow V$, its kernel pair (if it exists) is an equivalence relation. Equivalence relations of the form \say{kernel pair of some morphism $f$} are called \textit{effective}. If coequalizers exist (at least those of kernel pairs), an equivalence relation is effective if and only if it is the kernel pair of the coequalizer of its projections.

\definition
A regular category is said to be \emph{Barr exact} if all its equivalence relations are effective.
\enddefinition

A diagram of product $p_0, p_1 \colon X \longrightarrow W$ (if this operation is available in $\mathbf{C}$) is an equivalence relation: the morphisms $r$, $s$ and $h$ of Definition~\ref{defn:eq} are induced by the universal property. Equivalence relations are then monomorphisms $E \longrightarrow X$ to which the above morphisms restrict.

%In our analysis, we are considering categories of the form $L\KFr_\lf$. Among monomorphisms, we have injective p-morphisms (the injective p-morphisms are exactly the regular monomorphisms and all monomorphisms are regular in case $L$ extends \textsf{K4}, by Proposition~\ref{prop:injsurj} and Proposition~\ref{prop:injsurj2}); injective p-morphisms can be identified (up to isomorphism) with generated subframes of their codomain, by taking the $\Set$-image. Given $W \in L\KFr_\lf$, let $p_0, p_1 \colon X \longrightarrow W$ be a diagram of product in $L\KFr_\lf$. It is straightforward to verify that a generated subframe $E \subseteq X$ is an equivalence relation (with the appropriate restrictions of $p_0$ and $p_1$ as projections) if and only if the morphisms in Definition~\ref{defn:eq} restrict to $E$.

Describing equivalence relations in $L\KFr_\lf$ over some $W \in L\KFr_\lf$ appears to be quite complicated. However, some of them can be nicely constructed by selecting, for each pair $(a,b)$ of elements of $W$, a set of admissible pairs of elements of $W$. Let $W \times W$ denote the product in $\Set$, and let $\mathcal{P}(W \times W)$ denote the power-set.
\definition\label{defn:select}
A \emph{pair selection} for $W$ is a function $\mathcal{A} \colon W \times W \longrightarrow \mathcal{P}(W \times W)$, associating to each pair $(a,b) \in W \times W$ a set of pairs $\mathcal{A}_{(a,b)} \subseteq W \times W$ (called the set of $(a,b)$-\emph{admissible pairs}), satisfying the following properties:
\begin{enumerate}
    \item $a' \in a^* ~~~~~\implies~~~~~ (a',a') \in \mathcal{A}_{(a,a)}$;
    \item $(a',b') \in \mathcal{A}_{(a,b)} ~~~~~\implies~~~~~ (b',a') \in \mathcal{A}_{(b,a)}$;
    \item $(a',b') \in \mathcal{A}_{(a,b)} ~\&~ (b',c') \in \mathcal{A}_{(b,c)} ~~~~~\implies~~~~~ (a',c') \in \mathcal{A}_{(a,c)}$.
\end{enumerate}
\enddefinition
Given a pair of p-morphisms $g_0, g_1 \colon U \longrightarrow W$%, for each $u \in U$, we denote by $\pi(u)$ the pair $(g_0(u),g_1(u)) \in W \times W$. If
and a pair selection $\mathcal{A}$ for $W$,% has been fixed,
we can define the subset
$$U_{\mathcal{A}} \coloneqq \{u \in U\ \vert\ \forall y \forall y' (y \in u^* ~\&~ y' \in y^* \implies (g_0(y'),g_1(y')) \in \mathcal{A}_{(g_0(y),g_1(y))})\},$$
which is obviously a generated subframe of $U$.

\proposition\label{prop:rel}
Let $L$ be a normal modal logic with the finite model property. Let $p_0, p_1 \colon X \longrightarrow W$ be a diagram of product in $L\KFr_\lf$ and let $\mathcal{A}$ be a pair selection for $W$. Then $X_\mathcal{A}$ (with restricted projections) is an equivalence relation over $W$ in $L\KFr_\lf$.
\endproposition
\proof
The inclusion of $X_{\mathcal{A}}$ into $X$ is a monomorphism in $L\KFr_\lf$ ($L\KFr_\lf$ is closed under generated subframes); as we observed before, in order to prove the claim, it suffice to check that the morphisms $r$, $s$ and $h$ in Definition~\ref{defn:eq} restrict to $X_\mathcal{A}$.

For r) we need to show the following: for each $w \in W$, we have $r(w) \in X_\mathcal{A}$. For this purpose, take $y \in r(w)^* = r(w^*)$, say $y = r(a)$ for some $a \in w^*$, and $y' \in y^* = r(a)^* = r(a^*)$, say $y' = r(a')$ for some $a' \in a^*$. By 1. in Definition~\ref{defn:select}, we have that
$$(p_0(y'),p_1(y')) = (a',a') \in \mathcal{A}_{(a,a)} = \mathcal{A}_{(p_0(y),p_1(y))}.$$

Condition s) can be rephrased as follows: for each $x \in X$, if $x \in X_\mathcal{A}$, then $s(x) \in X_\mathcal{A}$. For this purpose, take $y \in s(x)^* = s(x^*)$, say $y = s(z)$ for some $z \in x^*$, and $y' \in y^* = s(z)^* = s(z^*)$, say $y' = s(z')$ for some $z' \in z^*$. Since $x \in X_\mathcal{A}$, we have that $(p_0(z'),p_1(z')) \in \mathcal{A}_{(p_0(z),p_1(z))}$. By 2. in Definition~\ref{defn:select}, we have that ($\sigma$ denotes the function $W \times W \longrightarrow W \times W$ switching the coordinates)
$$(p_0(y'),p_1(y')) = (p_1(z'),p_0(z')) \in \mathcal{A}_{(p_1(z),p_0(z))} = \mathcal{A}_{(p_0(y),p_1(y))}.$$

It remains to prove t). Let $h_0, h_1 \colon T \longrightarrow X$ be morphisms in $L\KFr_\lf$ such that $p_0 h_1 = p_1 h_0$ and let $h \colon T \longrightarrow X$ be the unique morphism such that $p_i h = p_i h_i$ for $i = 0, 1$. We need to show that, for each $t \in T$, if $h_0(t), h_1(t) \in X_\mathcal{A}$, then $h(t) \in X_\mathcal{A}$. Similarly to r) and s), using 3. in Definition~\ref{defn:select} we obtain our claim.
\endproof

We can use the previous result to obtain necessary conditions for $L\KFr_\lf$ to be Barr exact.% Let $g_0, g_1 \colon U \longrightarrow W$ be p-morphisms in $L\KFr_\lf$ and let $\mathcal{A}$ be a pair selection for $W$.% Denote by $f_\mathcal{A} \colon W \longrightarrow V_\mathcal{A}$ the coequalizer of the restrictions of $g_0$ and $g_1$ to the generated subframe $U_\mathcal{A}$.
\corollary\label{cor:rel}
Let $L$ be a normal modal logic with the finite model property. Assume that there is a pair of parallel p-morphisms $g_0, g_1 \colon U \longrightarrow W$ in $L\KFr_\lf$ and a pair selection $\mathcal{A}$ for $W$ such that the following conditions hold:
\begin{enumerate}
    \item $U \setminus U_\mathcal{A}$ is non-empty;
    \item for each $u \in U \setminus U_\mathcal{A}$, the pair $(g_0(u),g_1(u))$ belongs to the (set-theoretic) equivalence relation generated by $\{(g_0(u'),g_1(u'))\ \vert\ u' \in U_\mathcal{A}\}$.%$f_\mathcal{A} g_0 = f_\mathcal{A} g_1$.
\end{enumerate}
Then $L\KFr_\lf$ is not Barr-exact.
\endcorollary
\proof
As before, we denote by $p_0, p_1 \colon X \longrightarrow W$ the product in $L\KFr_\lf$ of $W$ with itself. Let $f \colon W \longrightarrow V$ be the coequalizer of $g_0$ and $g_1$. The kernel pair of $f$ is the equalizer of $fp_0$ and $fp_1$, hence, by Lemma~\ref{lem:eq}, it can be identified with a generated subframe $K \subseteq X$. The pair $g_0, g_1 \colon U \longrightarrow W$ determines a morphism $g \colon U \longrightarrow K \subseteq X$ such that $p_i g = g_i$, for $i = 0,1$, by the universal property of the kernel pair.

Consider now $X_\mathcal{A} \cap K$; it is an equivalence relation over $W$ in $L\KFr_\lf$, by Proposition~\ref{prop:rel} and by the fact that equivalence relations are closed under intersections. It is straightforward to verify that, for each $u \in U$, we have $u \in U_\mathcal{A}$ if and only if $g(u) \in X_\mathcal{A} \cap K$.

By 1. we have that $K \setminus (X_\mathcal{A} \cap K)$ is non-empty: there exists $u \in U \setminus U_\mathcal{A}$ and therefore $g(u) \in K \setminus (X_\mathcal{A} \cap K)$.

We now show that 2. implies that the coequalizer of $X_\mathcal{A} \cap K$ in $L\KFr_\lf$ is $f$: since the forgetful functor $L\KFr_\lf \longrightarrow \Set$ creates colimits, it suffices to prove that $\{(p_0(x),p_1(x))\ \vert\ x \in X_\mathcal{A} \cap K\}$ and $\{(g_0(u),g_1(u))\ \vert\ u \in U\}$ generate the same set-theoretic equivalence relation over $W$. If $x \in X_\mathcal{A} \cap K$, then $fp_0(x) = fp_1(x)$ (because $K$ is the kernel pair of $f$), hence $(p_0(x),p_1(x))$ belongs to the equivalence relation generated by $\{(g_0(u),g_1(u))\ \vert\ u \in U\}$. Conversely, if $u \in U$, then $(g_0(u),g_1(u))$ belongs to the equivalence relation generated by $\{(g_0(u'),g_1(u'))\ \vert\ u' \in U_\mathcal{A}\}$ (by 2.), and this set of generators is contained in $\{(p_0(x),p_1(x))\ \vert\ x \in X_\mathcal{A} \cap K\}$ because, for each $u' \in U_\mathcal{A}$, we have $g(u') \in X_\mathcal{A} \cap K$ and $(g_0(u'),g_1(u')) = (p_0g(u'),p_1g(u'))$.

Therefore, $X_\mathcal{A} \cap K$ is an equivalence relation in $L\KFr_\lf$ over $W$ that is not effective (it is not isomorphic to the kernel pair $K$ of its coequalizer $f$).
\endproof

Applying the previous result, we can prove that Barr exactness excludes a lot of configurations.
\proposition\label{prop:necexact}
Let $L$ be a normal modal logic with the finite model property. If $L\KFr_\lf$ is Barr exact, then it cannot contain any of the following Kripke frames: %$[4]$, $\fork{2}$, $C_3$, $(C_2)^+$ and $(C_3)^-$.
\[
    \begin{tikzpicture}
    \filldraw (0,0.6) circle (0.1cm);
    \filldraw (0,0.2) circle (0.1cm);
    \filldraw (0,-0.2) circle (0.1cm);
    \filldraw (0,-0.6) circle (0.1cm);
    \draw (0,-0.5) -- (0,-0.3) (0,-0.1) -- (0,0.1) (0,0.3) -- (0,0.5);
    \begin{scope}[shift={(2.5,0)}]
    \filldraw (-0.2,0.2) circle (0.1cm);
    \filldraw (0.2,0.2) circle (0.1cm);
    \filldraw (0,-0.2) circle (0.1cm);
    \draw (0,-0.2) -- (-0.2,0.2) (0,-0.2) -- (0.2,0.2);
    \begin{scope}[shift={(2.5,0)}]
    \draw (0,0) circle (0.5cm);
    \filldraw (-0.2,0.1) circle (0.1cm);
    \filldraw (0.2,0.1) circle (0.1cm);
    \filldraw (0,-0.2) circle (0.1cm);
    \begin{scope}[shift={(2.5,0)}]
    \filldraw (0,-0.5) circle (0.1cm);
    \draw (0,0.2) circle (0.4cm);
    \filldraw (-0.15,0.2) circle (0.1cm);
    \filldraw (0.15,0.2) circle (0.1cm);
    \draw (0,-0.4) -- (0,-0.2);
    \begin{scope}[shift={(2.5,0)}]
    \draw (0,-0.2) circle (0.5cm);
    \filldraw (-0.2,-0.1) circle (0.1cm);
    \filldraw (0.2,-0.1) circle (0.1cm);
    \filldraw (0,-0.4) circle (0.1cm);
    \filldraw (0,0.6) circle (0.1cm);
    \draw (0,0.5) -- (0,0.3);
    \end{scope}
    \end{scope}
    \end{scope}
    \end{scope}
    \end{tikzpicture}
\]
\endproposition
\proof
We apply Corollary~\ref{cor:rel}. As an example, assume that $[4] \in L\KFr_\lf$. Let $g_0, g_1 \colon U \longrightarrow W$ be the pair of parallel p-morphisms in $L\KFr_\lf$ depicted below
\begin{comment}
\[
    \begin{tikzpicture}
    \filldraw[red] (0,0.6) circle (0.1cm);
    \node[anchor=east] at (0,0.6) {$(1,1)$};
    \filldraw[red] (0,0.2) circle (0.1cm);
    \node[anchor=east] at (0,0.2) {$(1,2)$};
    \filldraw[red] (0,-0.2) circle (0.1cm);
    \node[anchor=east] at (0,-0.2) {$(2,2)$};
    \filldraw (0,-0.6) circle (0.1cm);
    \node[anchor=east] at (0,-0.6) {$(3,3)$};
    \draw (0,-0.5) -- (0,-0.3) (0,-0.1) -- (0,0.1) (0,0.3) -- (0,0.5);
    \begin{scope}[shift={(2,0)}]
    \draw[->] (-0.8,0) -- (0.8,0) node [midway, above=1pt, fill=white] {};
    \begin{scope}[shift={(2,0)}]
    \filldraw (0,0.4) circle (0.1cm);
    \node[anchor=west] at (0,0.4) {$1$};
    \filldraw (0,0) circle (0.1cm);
    \node[anchor=west] at (0,0) {$2$};
    \filldraw (0,-0.4) circle (0.1cm);
    \node[anchor=west] at (0,-0.4) {$3$};
    \draw (0,-0.3) -- (0,-0.1) (0,0.1) -- (0,0.3);
    \end{scope}
    \end{scope}
    \end{tikzpicture}
\]
\end{comment}
\begin{center}
\begin{tikzpicture}[
  vertex/.style = {circle, draw,
  minimum size=26, inner sep=0,
  fill=white},
  vertex1/.style = {vertex, fill=red!20!white},
  vertex2/.style = {vertex, fill=orange!30!white},
  vertex3/.style = {vertex, fill=blue!30!white},
  vertex4/.style = {vertex, fill=teal!30!white},
]

  \draw[thick]
    (-2,3) node[vertex1] (1) {$1,1$}
    (-2,1) node[vertex1] (2) {$1,2$}
    (-2,-1) node[vertex1] (3) {$2,2$}
    (-2,-3) node[vertex] (4) {$3,3$}
    (4) -- (3)
    (3) -- (2)
    (2) -- (1);

   \draw[thick]
    (2,2) node[vertex] (a) {$1$}
    (2,0) node[vertex] (b) {$2$}
    (2,-2) node[vertex] (c) {$3$}
    (c) -- (b)
    (b) -- (a);

   \draw[->] (-1,0) -- (1,0);
\end{tikzpicture}
\end{center}
where the red dots represent the generated subframe $U_\mathcal{A}$ determined by the pair selection $\mathcal{A}$ for $W$ such that: $\mathcal{A}_{(a,b)}$ contains all the pairs in $W \times W$, if both $a$ and $b$ are different from $3$; $\mathcal{A}_{(a,b)}$ coincides with the diagonal $\{(1,1),(2,2),(3,3)\}$, otherwise. Conditions 1. and 2. of Corollary~\ref{cor:rel} are easily verified.

The pictures below illustrate the pairs of p-morphisms and the corresponding pair selections that can be chosen in order to prove the claim for the remaining cases, using Corollary~\ref{cor:rel}.
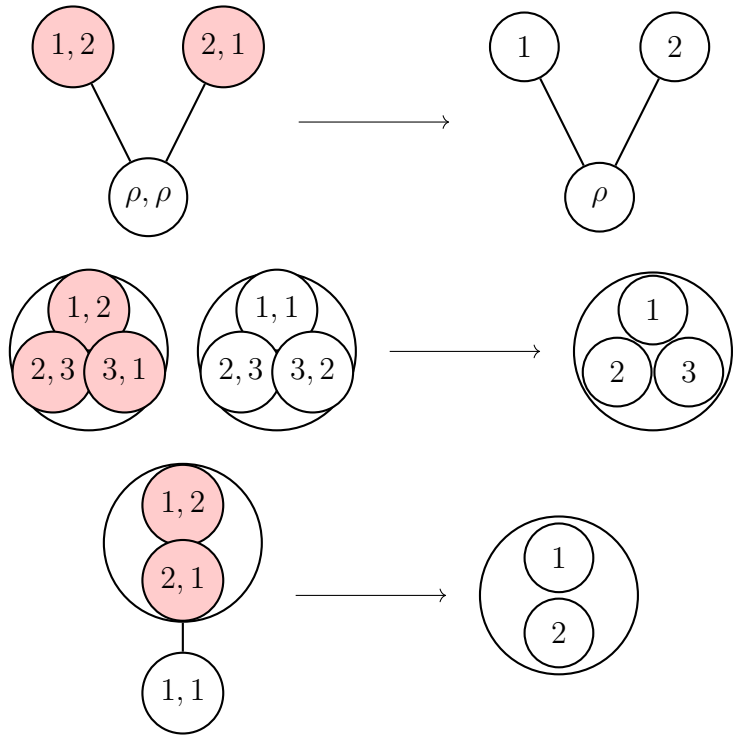
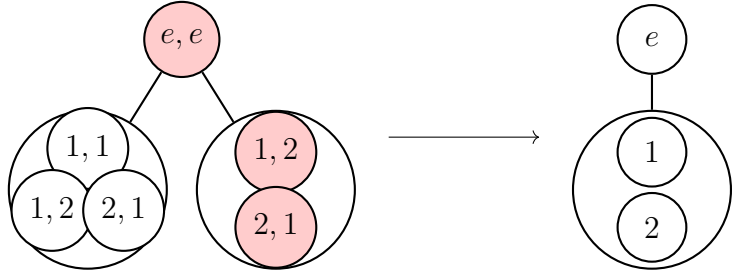
\begin{center}
\begin{tikzpicture}[
  vertex/.style = {circle, draw,
  minimum size=26, inner sep=0,
  fill=white},
  vertex1/.style = {vertex, fill=red!20!white},
  vertex2/.style = {vertex, fill=orange!30!white},
  vertex3/.style = {vertex, fill=blue!30!white},
  vertex4/.style = {vertex, fill=teal!30!white},
]

  \draw[thick]
    (-4,1) node[vertex1] (1) {$1,2$}
    (-2,1) node[vertex1] (2) {$2,1$}
    (-3,-1) node[vertex] (r) {$\rho,\rho$}
    (r) -- (1)
    (r) -- (2);

  \draw[thick]
    (2,1) node[vertex] (1') {$1$}
    (4,1) node[vertex] (2') {$2$}
    (3,-1) node[vertex] (r') {$\rho$}
    (r') -- (1')
    (r') -- (2');

   \draw[->] (-1,0) -- (1,0);
\end{tikzpicture}
\end{center}

\begin{center}
\begin{tikzpicture}[
  vertex/.style = {circle, draw,
  minimum size=26, inner sep=0,
  fill=white},
  cluster/.style = {circle, draw,
  minimum size=22mm, inner sep=0,
  fill=white},
  vertex1/.style = {vertex, fill=red!20!white},
  vertex2/.style = {vertex, fill=orange!30!white},
  vertex3/.style = {vertex, fill=blue!30!white},
  vertex4/.style = {vertex, fill=teal!30!white},
]

  \draw[thick]
    (-3.75,0) node[cluster] (A) {}
    ([shift=({270:0.55})]A) node[vertex1] {$1,2$}
    ([shift=({150:0.55})]A) node[vertex1] {$2,3$}
    ([shift=({30:0.55})]A) node[vertex1] {$3,1$}
    (-1.25,0) node[cluster] (B) {}
    ([shift=({270:0.55})]B) node[vertex] {$1,1$}
    ([shift=({150:0.55})]B) node[vertex] {$2,3$}
    ([shift=({30:0.55})]B) node[vertex] {$3,2$};

  \draw[thick]
    (3.75,0) node[cluster] (C) {}
    ([shift=({270:0.55})]C) node[vertex] {$1$}
    ([shift=({150:0.55})]C) node[vertex] {$2$}
    ([shift=({30:0.55})]C) node[vertex] {$3$};

   \draw[->] (0.25,0) -- (2.25,0);
\end{tikzpicture}
\end{center}

\begin{center}
\begin{tikzpicture}[
  vertex/.style = {circle, draw,
  minimum size=26, inner sep=0,
  fill=white},
  cluster/.style = {circle, draw,
  minimum size=22mm, inner sep=0,
  fill=white},
  vertex1/.style = {vertex, fill=red!20!white},
  vertex2/.style = {vertex, fill=orange!30!white},
  vertex3/.style = {vertex, fill=blue!30!white},
  vertex4/.style = {vertex, fill=teal!30!white},
]

  \draw[thick]
    (-2.5,0.7) node[cluster] (A) {}
    ([shift=({180:0.5})]A) node[vertex1] {$1,2$}
    ([shift=({0:0.5})]A) node[vertex1] {$2,1$}
    (-2.5,-1.3) node[vertex] (B) {$1,1$}
    (A) -- (B);

  \draw[thick]
    (2.5,0) node[cluster] (C) {}
    ([shift=({180:0.5})]C) node[vertex] {$1$}
    ([shift=({0:0.5})]C) node[vertex] {$2$};

  \draw[->] (-1,0) -- (1,0);
\end{tikzpicture}
\end{center}

\begin{center}
\begin{tikzpicture}[
  vertex/.style = {circle, draw,
  minimum size=26, inner sep=0,
  fill=white},
  cluster/.style = {circle, draw,
  minimum size=22mm, inner sep=0,
  fill=white},
  vertex1/.style = {vertex, fill=red!20!white},
  vertex2/.style = {vertex, fill=orange!30!white},
  vertex3/.style = {vertex, fill=blue!30!white},
  vertex4/.style = {vertex, fill=teal!30!white},
]

  \draw[thick]
    (-3.125,1.3) node[vertex1] (ee) {$\varepsilon,\varepsilon$}
    (-4.375,-0.7) node[cluster] (A) {}
    ([shift=({270:0.55})]A) node[vertex] {$1,1$}
    ([shift=({150:0.55})]A) node[vertex] {$1,2$}
    ([shift=({30:0.55})]A) node[vertex] {$2,1$}
    (-1.875,-0.7) node[cluster] (B) {}
    ([shift=({180:0.5})]B) node[vertex1] {$1,2$}
    ([shift=({0:0.5})]B) node[vertex1] {$2,1$}
    (ee) -- (A)
    (ee) -- (B);

  \draw[thick]
    (3.125,1.3) node[vertex] (e) {$\varepsilon$}
    (3.125,-0.7) node[cluster] (C) {}
    ([shift=({180:0.5})]C) node[vertex] {$1$}
    ([shift=({0:0.5})]C) node[vertex] {$2$}
    (e) -- (C);

  \draw[->] (-0.375,0) -- (1.625,0);
\end{tikzpicture}
\end{center}
\endproof

Recall the list provided by Theorem~\ref{teo:regchar} of the logics $L \in \NExt\textsf{S4}$ with the finite model property such that $L\KFr_\lf$ is a regular category (regularity is a necessary condition for Barr exactness).
\theorem\label{teo:exact}
Let $L$ be a normal modal logic with the finite model property and assume that $L \in \NExt \textup{\textsf{S4}}$. If $L\KFr_\lf$ is a Barr exact category, then $L$ is the inconsistent logic, or $L\KFr_\lf$ is one of the following:
\begin{enumerate}
    \item[(i)] $\{W \in \textup{\textsf{S4}}\KFr_\lf\ \vert\ \depth(W) \leq 2,\ \width(W) \leq 1,\ \exter(W) \leq 1,\ \inter(W) \leq 1\}$;
    \item[(ii)] $\{W \in \textup{\textsf{S4}}\KFr_\lf\ \vert\ \depth(W) \leq 2,\ \width(W) \leq 1,\ \exter(W) \leq 1,\ \inter(W) \leq 2\}$;
    \item[(iii)] $\{W \in \textup{\textsf{S4}}\KFr_\lf\ \vert\ \depth(W) \leq 1,\ \exter(W) \leq 1\}$;
    \item[(iv)] $\{W \in \textup{\textsf{S4}}\KFr_\lf\ \vert\ \depth(W) \leq 1,\ \exter(W) \leq 2\}$.
\end{enumerate}
\endtheorem
\proof
From the logics in 4. of Theorem~\ref{teo:regchar}, exclude those containing any of the Kripke frames of Proposition~\ref{prop:necexact}.
\endproof

As we did for regularity, it is natural to ask whether the classification of Theorem~\ref{teo:exact} is complete. This is answered in the next section: in all cases $L\KFr_\lf$ turns out to be equivalent to the category of actions of some finite monoid, and an instance of a general adjunction between a category of presheaves and $\KFr$.

\section{An adjunction}\label{sec:presheaves}
When $L\KFr_\lf \simeq (\Pro L\MA_\fin)^\op$ is Barr exact, we can say more. As we said in the introduction, $L\KFr_\lf$ has all colimits, inherited by the full embedding $L\KFr_\lf \longrightarrow \KFr$, and the forgetful functor $L\KFr_\lf \longrightarrow \Set$ creates them (see Proposition~\ref{prop:colim}). In particular, coproducts in $L\KFr_\lf$ are computed as disjoint unions. Namely, given a discrete diagram $(W_i \mid i \in I)$ in $L\KFr_\lf$, the coprojections $\iota_j \colon W_j \longrightarrow \coprod_{i \in I} W_i$ are given by the inclusions of the $W_j$'s into their disjoint union (equipped with the disjoint union of the binary relations).
\definition
A category with coproducts is said to be \emph{infinitary extensive} if
\begin{enumerate}
    \item pullbacks of coprojections along arbitrary morphisms exist;
    \item coproducts are \emph{disjoint} (coprojections are monomorphisms and the pullback of two different coprojections is the initial object);
    \item coproducts are stable under pullback.
\end{enumerate}
\enddefinition

\proposition\label{prop:ext}
Let $L$ be a normal modal logic with the finite model property. The category $L\KFr_\lf$ is infinitary extensive.
\endproposition
\proof
Each coprojection $\iota_j \colon W_j \longrightarrow \coprod_{i \in I} W_i$ is a monomorphism in $L\KFr_\lf$, being an injective p-morphism. Moreover, it is straightforward to see that the following
\[\begin{tikzcd}
	\emptyset & {W_k} \\
	{W_j} & {\coprod_{i \in I} W_i}
	\arrow[from=1-1, to=1-2]
	\arrow[from=1-1, to=2-1]
	\arrow["{\iota_k}", from=1-2, to=2-2]
	\arrow["{\iota_j}"', from=2-1, to=2-2]
\end{tikzcd}\]
is a pullback diagram in $L\KFr_\lf$ if $j \neq k \in I$, and that $\emptyset$ is the initial object of $L\KFr_\lf$.

The $L\KFr_\lf$-pullback of a coprojection $\iota_j \colon W_j \longrightarrow \coprod_{i \in I} W_i$ along any morphism $f \colon V \longrightarrow \coprod_{i \in I} W_i$ is given by the the following diagram
\[\begin{tikzcd}
	{f^{-1}(W_j)} & V \\
	{W_j} & {\coprod_{i \in I} W_i}
	\arrow[from=1-1, to=1-2]
	\arrow[from=1-1, to=2-1]
	\arrow["f", from=1-2, to=2-2]
	\arrow["{\iota_j}"', from=2-1, to=2-2]
\end{tikzcd}\]
where $f^{-1}(W_j) \longrightarrow V$ is the inclusion and $f^{-1}(W_j) \longrightarrow W_j$ is the restriction of $f$. $(f^{-1}(W_j) \longrightarrow V\ \vert\ i \in I)$ is a coproduct diagram, since $V$ is the disjoint union of the $f^{-1}(W_j)$'s.
\endproof

Observe also that $L\KFr_\fin$ (the class of the finite Kripke frames belonging to $L\KFr$) is a \emph{generating family} for $L\KFr_\lf$\footnote{A family $\mathbf{G}$ of objects of a category $\mathbf{C}$ is said to be \emph{generating} if two parallel morphisms coincide whenever they are equalized by all (composable) morphisms having domain in $\mathbf{G}$.}, since $L\KFr_\lf$ is the Ind-completion of $L\KFr_\fin$ via the full inclusion $L\KFr_\fin \subseteq L\KFr_\lf$, as stated in Theorem~\ref{teo:ind}. Moreover, $L\KFr_\fin$ is \emph{essentially small} (it is equivalent to a small category). Using Giraud's characterization (see \cite{MR2522659}) of \emph{Grothendieck toposes} (sheaves for a certain site) as Barr exact and infinitary extensive categories with a small set of generators, we obtain the following.
\theorem
Let $L$ be a normal modal logic with the finite model property. The category $L\KFr_\lf$ is Barr exact if and only if it is a Grothendieck topos.
\endtheorem
\proof
If follows from Proposition~\ref{prop:ext} and the previous observations.
\endproof

We want to prove that all the $L$ appearing in the list of Theorem~\ref{teo:exact} are such that $L\KFr_\lf$ is Barr exact, by explicitly describing $L\KFr_\lf$ as a Grothendieck topos. Among toposes, there are categories of presheaves.

There is a canonical way to see presheaves as Kripke frames. Given a small category $\mathbf{C}$ and a presheaf $F \in \widehat{\mathbf{C}}$\footnote{The category of presheaves $\mathbf{C}^\op \longrightarrow \Set$ and natural transformations between them.}, we can consider the slice category $\mathbf{C}/F$ (called \emph{diagram of} $F$, or \emph{category of the elements of} $F$). The objects of $\mathbf{C}/F$ are natural transformations $x \colon \mathbf{C}[-,c] \longrightarrow F$ from representable presheaves to $F$\footnote{They form a set isomorphic to the disjoint union of the $F(c)$'s, with $c$ in $\mathbf{C}$ (by Yoneda's lemma).}; morphisms $y \longrightarrow x$ in $\mathbf{C}/F$ are given by arrows $h \colon d \longrightarrow c$ in $\mathbf{C}$ such that
$$y = x \mathbf{C}[-,h],$$
where $\mathbf{C}[-,h] \colon \mathbf{C}[-,d] \longrightarrow \mathbf{C}[-,c]$ is the natural transformation acting by composition with $h$. Seeing $\mathbf{C}$ as a full subcategory of $\widehat{\mathbf{C}}$, we will make no distinction between $c \in \mathbf{C}$ and the representable presheaf $\mathbf{C}[-,c]$, as well as between $h \colon d \longrightarrow c$ in $\mathbf{C}$ and the aforementioned natural transformation.

We can then consider the Kripke frame $\mathcal{K}(F)$, with the objects of $\mathbf{C}/F$ as underlying set, obtained by setting
$$x \prec y \iff y = x h \text{ for some } h \colon d \longrightarrow c \text{ in } \mathbf{C},$$
for $x \colon c \longrightarrow F$ and $y \colon d \longrightarrow F$ in $\mathbf{C}/F$. Observe that $\mathcal{K}(F)$ is a preorder (i.e.\ $\mathcal{K}(F) \in \textsf{S4}\KFr$), which is known as \emph{preordered reflection} of the category $\mathbf{C}/F$ (in \cite{Ghilardi1992QuantifiedEO}\footnote{Used to provide a complete generalized Kripke semantics for quantified extensions of canonical propositional intermediate logics.} it is called \emph{frame representation}). Observe that, for $x \colon c \longrightarrow F$ in $\mathcal{K}(F)$,
$$\up{x} = \{x h\ \vert\ h \colon d \longrightarrow c \text{ in } \mathbf{C}\}.$$

A natural transformation $\alpha \colon F \longrightarrow G$ induces a functor $\alpha_* \colon \mathbf{C}/F \longrightarrow \mathbf{C}/G$ acting by composition. Seen as a function $\mathcal{K}(F) \longrightarrow \mathcal{K}(G)$, $\alpha_*$ is a p-morphism: for $x \colon c \longrightarrow F$ in $\mathcal{K}(F)$,
\begin{align*}
\alpha_*(\up{x}) &= \alpha_*(\{x h\ \vert\ h \colon d \longrightarrow c \text{ in } \mathbf{C}\}) = \\
&= \{\alpha x h\ \vert\ h \colon d \longrightarrow c \text{ in } \mathbf{C}\} = \up{(\alpha x)} = \up{\alpha_*(x)}.
\end{align*}

Putting everything together,
\lemma\label{lem:functor}
$\mathcal{K}$ defines a faithful functor $\widehat{\mathbf{C}} \longrightarrow \KFr$.
\endlemma
\proof
The functoriality part follows from the previous observations. Faithfullness follows from the fact that representable presheaves in $\widehat{\mathbf{C}}$ form a generating family (by Yoneda's lemma), hence the inequality of $\beta, \gamma \colon F \longrightarrow G$ in $\widehat{\mathbf{C}}$ is testified as $\mathcal{K}(\beta)(x) = \beta x \neq \gamma x = \mathcal{K}(\gamma)(x)$ for some $x \colon c \longrightarrow F$ in $\mathbf{C}/F = \mathcal{K}(F)$.
\endproof

It is a well known fact in category theory (\cite{CWM} for reference) that a morphism $\alpha \colon F \longrightarrow G$ in $\widehat{\mathbf{C}}$ is a monomorphism (resp.\ epimorphism) if and only if all its components are injective (resp.\ surjective) functions.
\lemma\label{lem:monoinj-episurj}
A morphism $\alpha \colon F \longrightarrow G$ in $\widehat{\mathbf{C}}$ is a monomorphism (resp.\ an epimorphism) if and only if the p-morphism $\mathcal{K}(\alpha) \colon \mathcal{K}(F) \longrightarrow \mathcal{K}(G)$ is injective (resp.\ surjective).
\endlemma
\proof
If $\alpha$ is a monomorphism, then $\mathcal{K}(\alpha) = \alpha_*$ is injective. Vice versa, assume that $\mathcal{K}(\alpha)$ is injective. Given $\beta, \gamma \colon H \longrightarrow F$ in $\widehat{\mathbf{C}}$ such that $\beta \neq \gamma$ and $\alpha \beta = \alpha \gamma$, fixing $x \colon c \longrightarrow H$ testifying $\beta \neq \gamma$ (i.e.\ $\beta x \neq \gamma x$ in $\mathcal{K}(F)$), we have that $\alpha \beta x = \mathcal{K}(\alpha)(\beta x) \neq \mathcal{K}(\alpha)(\gamma x) = \alpha \gamma x$ (against $\alpha \beta = \alpha \gamma$).

If $\alpha$ is an epimorphism, then $\mathcal{K}(\alpha) = \alpha_*$ is surjective. Indeed, given $y \colon d \longrightarrow G$ in $\mathcal{K}(G)$, using projectivity of the representables in $\widehat{\mathbf{C}}$ (by Yoneda's lemma), we can find $x \colon c \longrightarrow F$ in $\mathcal{K}(F)$ such that $\mathcal{K}(\alpha)(x) = \alpha x = y$. Vice versa, assume that $\mathcal{K}(\alpha)$ is surjective. Given $\beta, \gamma \colon G \longrightarrow H$ such that $\beta \neq \gamma$ and $\beta \alpha = \gamma \alpha$, fixing $y \colon d \longrightarrow G$ in $\mathcal{K}(G)$ testifying $\beta \neq \gamma$ (i.e.\ $\beta y \neq \gamma y$), we can find some $x \colon c \longrightarrow F$ in $\mathcal{K}(F)$ such that $\alpha x = \mathcal{K}(\alpha)(x) = y$. As a consequence, $\beta y = \beta \alpha x = \gamma \alpha x = \gamma y$ (against $\beta y \neq \gamma y$).
\endproof

Going in the other direction, we can associate a presheaf in $\widehat{\mathbf{C}}$ to any Kripke frame $W$. For $c$ in $\mathbf{C}$, consider the set $\mathcal{F}(W)(c) \coloneqq \KFr[\mathcal{K}(c),W]$ of the p-morphisms from the preorder $\mathcal{K}(c)$ to $W$\footnote{Recall that we are identifying $c \in \mathbf{C}$ with the presheaf $\mathbf{C}[-,c]$, so that $\mathcal{K}(c)$ denotes $\mathcal{K}(\mathbf{C}[-,c])$ and it is isomorphic to the set of arrows $h \colon d \longrightarrow c$ in $\mathbf{C}$ having $c$ as codomain. With such identifications, for $h \colon d \longrightarrow c$ in $\mathbf{C}$, the p-morphism $\mathcal{K}(h) \colon \mathcal{K}(d) \longrightarrow \mathcal{K}(c)$ acts by composition.}. For $h \colon d \longrightarrow c$ in $\mathbf{C}$, define the function $\mathcal{F}(W)(h) \colon \mathcal{F}(W)(c) \longrightarrow \mathcal{F}(W)(d)$ sending the p-morphism $g \colon \mathcal{K}(c) \longrightarrow W$ to the p-morphism $g \mathcal{K}(h) \colon \mathcal{K}(d) \longrightarrow W$. It is straightforward to verify that $\mathcal{F}(W)$ defines a presheaf in $\widehat{\mathbf{C}}$.

To a p-morphism $f \colon W \longrightarrow V$, we associate the natural transformation $\mathcal{F}(f)$, whose component $\mathcal{F}(f)_c \colon \mathcal{F}(W)(c) \longrightarrow \mathcal{F}(V)(c)$ sends the p-morphism $g \colon \mathcal{K}(c) \longrightarrow W$ to the p-morphism $fg \colon \mathcal{K}(c) \longrightarrow V$.

\lemma
$\mathcal{F}$ defines a functor $\KFr \longrightarrow \widehat{\mathbf{C}}$.
\endlemma
\proof
Easy to check.
\endproof

For each Kripke frame $W$, we can define a function $\varepsilon_W \colon \mathcal{K}(\mathcal{F}(W)) \longrightarrow W$ as follows: via the natural isomorphism $\widehat{\mathbf{C}}[c,\mathcal{F}(W)] \cong \mathcal{F}(W)(c) = \KFr[\mathcal{K}(c),W]$,
\[\begin{tikzcd}[row sep=tiny]
	{\coprod_{c \in \mathbf{C}} \KFr[\mathcal{K}(c),W]} & {W} \\
	{g \colon \mathcal{K}(c) \longrightarrow W} & {g(\id_c)}
	\arrow["{\varepsilon_W}", from=1-1, to=1-2]
	\arrow[maps to, from=2-1, to=2-2]
\end{tikzcd}\]
The function $\varepsilon_W$ is a p-morphism:
\begin{align*}
\varepsilon_W(\up{g}) &= \varepsilon_W(\{g\mathcal{K}(h)\ \vert\ h \colon d \longrightarrow c \text{ in } \mathbf{C}\}) = \{g\mathcal{K}(h)(\id_{d})\ \vert\ h \colon d \longrightarrow c \text{ in } \mathbf{C}\} = \\
&= \{g(h)\ \vert\ h \colon d \longrightarrow c \text{ in } \mathbf{C}\} = g(\up{\id_c}) = \up{g(\id_c)} = \up{\varepsilon_W(g)}
\end{align*}
The $\varepsilon_W$'s are the components of a natural transformation $\varepsilon \colon \mathcal{K} \mathcal{F} \longrightarrow \id_{\KFr}$.

For each presheaf $F \in \widehat{\mathbf{C}}$ and $c \in \mathbf{C}$, we can define a function $(\eta_F)_c \colon F(c) \longrightarrow \mathcal{F}(\mathcal{K}(F))(c)$ as follows: via the natural isomorphism $\widehat{\mathbf{C}}[c,F] \cong F(c)$,
\[\begin{tikzcd}[row sep=tiny]
	{\widehat{\mathbf{C}}}[c,F] & {\KFr[\mathcal{K}(c),\mathcal{K}(F)]} \\
	x & {\mathcal{K}(x)}
	\arrow["{(\eta_F)_c}", from=1-1, to=1-2]
	\arrow[maps to, from=2-1, to=2-2]
\end{tikzcd}\]
The functions $(\eta_F)_c$ are the component of a natural transformation $\eta_F \colon F \longrightarrow \mathcal{F}(\mathcal{K}(F))$: for $h \colon d \longrightarrow c$ in $\mathbf{C}$,
\begin{align*}
\mathcal{F}(\mathcal{K}(F))(h)(\eta_F)_c(x) &= \mathcal{F}(\mathcal{K}(F))(h)(\mathcal{K}(x)) = \\
&= \mathcal{K}(x) \mathcal{K}(h) = \mathcal{K}(xh) = (\eta_F)_d(xa) = (\eta_F)_d F(h)(x)
\end{align*}
The $\eta_F$'s are the components of a natural transformation $\eta \colon \id_{\widehat{\mathbf{C}}} \longrightarrow \mathcal{F} \mathcal{K}$.
\proposition
$\mathcal{K}$ is left adjoint to $\mathcal{F}$, with $\eta$ and $\varepsilon$ respectively unit and counit of the adjunction.
\endproposition
\proof
It is enough to verify the triangle identities that characterize an adjoint pair of functors (see \cite{CWM}).
\endproof

The construction above allows us to characterize fullness for $\mathcal{K}$.
\lemma\label{lem:counit}
$\mathcal{K}$ is full if and only if the p-morphism $\varepsilon_{\mathcal{K}(F)}$ is an isomorphism for all $F \in \widehat{\mathbf{C}}$.
\endlemma
\proof
By a well known result in category theory (see \cite{CWM}), a left adjoint is full and faithful (and we observed in Lemma~\ref{lem:functor} that $\mathcal{K}$ is always faithful) if and only if the unit of the adjunction is a natural isomorphism (all the components $\eta_F$, with $F \in \widehat{\mathbf{C}}$, are isomorphisms).

By Lemma~\ref{lem:monoinj-episurj}, we have that $\eta_F \colon F \longrightarrow \mathcal{F}(\mathcal{K}(F))$ in $\widehat{\mathbf{C}}$ is an isomorphism if and only if its image $\mathcal{K}(\eta_F) \colon \mathcal{K}(F) \longrightarrow \mathcal{K}(\mathcal{F}(\mathcal{K}(F)))$ in $\KFr$ is so.

The claim now follows from one of the triangle identities satisfied by the unit and the counit of an adjoint pair, namely that
\[\begin{tikzcd}
	{\mathcal{K}(\mathcal{F}(\mathcal{K}(F)))} & {\mathcal{K}(F)} \\
	{\mathcal{K}(F)}
	\arrow["{\varepsilon_{\mathcal{K}(F)}}", from=1-1, to=1-2]
	\arrow["{\mathcal{K}(\eta_F)}", from=2-1, to=1-1]
	\arrow[Rightarrow, no head, from=2-1, to=1-2]
\end{tikzcd}\]
commutes.
\endproof

The previous triangle identity tells us that $\varepsilon_{\mathcal{K}(F)}$ is always a surjective p-morphism, for all $F \in \widehat{\mathbf{C}}$. A sufficient condition for the fullness of $\mathcal{K} \colon \widehat{\mathbf{C}} \longrightarrow \KFr$ is then that $\varepsilon_W$ is injective, for all Kripke frames $W$. Such a sufficient condition allows for an explicit characterization of the class of Kripke frames obtained as the image of $\mathcal{K}$, provided that the following additional condition is taken into account: in order to guarantee the locally finiteness of each $\mathcal{K}(F)$, from now on, we will assume that all the slices $\mathbf{C}/c$, with $c \in \mathbf{C}$, are finite (recall that, for $F \in \widehat{\mathbf{C}}$ and $x \colon c \longrightarrow F \in \mathcal{K}(F)$, we have $x^* = \up{x} = \{x h\ \vert\ h \colon d \longrightarrow c \text{ in } \mathbf{C}\}$).

Let $L_\mathbf{C}$ be the logic of the class of Kripke frames $\{\mathcal{K}(c)\ \vert\ c \in \mathbf{C}\}$, i.e.\ the set of modal formulas valid in every $\mathcal{K}(c)$. Then $L_\mathbf{C} \in \NExt \textsf{S4}$ and it has the finite model property (each $\mathcal{K}(c)$ is a finite preorder). A characterization of modal definability, similar to that for finite transitive Kripke frames in \cite[Theorem 3.21]{MR1837791}, holds for locally finite frames: a class of locally finite transitive Kripke frames is of the form $L\KFr_\lf$, for some $L \in \NExt \textsf{K4}$, if and only if it is closed under generated subframes, p-morphic images and arbitrary disjoint unions. This means that $L_\mathbf{C}\KFr_\lf$ is obtained by closing $\{\mathcal{K}(c)\ \vert\ c \in \mathbf{C}\}$ under the operations above. We write $V \subred W$ (to be read as $V$ \emph{subreduces to} $W$) to denote the fact that there exists a generated subframe $G \subseteq V$ and a surjective p-morphism $G \longrightarrow W$. We then have that
$$L_\mathbf{C}\KFr_\lf = \{W \in \KFr_\lf\ \vert\ (\forall w \in W)(\exists c \in \mathbf{C} \text{ s.t.\ } \mathcal{K}(c) \subred w^*)\}.$$
Indeed, the right-hand side is easily seen to be the smallest class of locally finite Kripke frames containing $\{\mathcal{K}(c)\ \vert\ c \in \mathbf{C}\}$ and satisfying the closure properties previously mentioned.
\theorem\label{teo:suff}
If $\varepsilon_W \colon \mathcal{K}(\mathcal{F}(W)) \longrightarrow W$ is an injective p-morphism for all $W \in \KFr$, then $\mathcal{K}$ restricts to an equivalence
$$\widehat{\mathbf{C}} \simeq L_{\mathbf{C}}\KFr_\lf.$$
%of categories between $\widehat{\mathbf{C}}$ and $L_{\mathbf{C}}\KFr_\lf$.%, where $L_{\mathbf{C}}$ is the logic of the class of finite Kripke frames $\{\mathcal{K}(c)\ \vert\ c \in \mathbf{C}\}$.
\endtheorem
\proof
The fullness of $\mathcal{K}$ follows from Lemma~\ref{lem:counit} and from the fact that $\varepsilon_{\mathcal{K}(F)}$ is always surjective, as we observed before. We need to prove that the image of $\mathcal{K}$ (closed under isomorphisms) is the class of Kripke frames $L_{\mathbf{C}}\KFr_\lf$.

If $F \in \widehat{\mathbf{C}}$, then $\mathcal{K}(F) \in L_\mathbf{C}\KFr_\lf$. Indeed, for each $x \colon c \longrightarrow F$ in $\mathcal{K}(F)$, the generated subframe $x^* = \{x h\ \vert\ h \colon d \longrightarrow c \text{ in } \mathbf{C}\}$ is the image of the p-morphism $\mathcal{K}(x) \colon \mathcal{K}(c) \longrightarrow \mathcal{K}(F)$.
%Consider $F$ in $\widehat{\mathbf{C}}$ and a propositional modal formula $\varphi \in L_{\mathbf{C}}$. If $\mathcal{K}(F) \nvDash \varphi$, then there exists $x \colon c \longrightarrow F$ in $\mathcal{K}(F)$ such that $x^* \nvDash \varphi$. However, this is in contradiction with the fact that $x^* = \{x h\ \vert\ h \colon d \longrightarrow c \text{ in } \mathbf{C}\}$ is the p-morphic image of $\mathcal{K}(c)$ via $\mathcal{K}(x) \colon \mathcal{K}(c) \longrightarrow \mathcal{K}(F)$ and that the logic is preserved by p-morphic images (see \cite{misha}).

Conversely, consider $W \in L_{\mathbf{C}}\KFr_\lf$; we prove that $\varepsilon_W \colon \mathcal{K}(\mathcal{F}(W)) \longrightarrow W$ is an isomorphism (it is sufficient to prove that it is surjective, since it is injective by hypothesis). Let us take $w \in W$; there exists some $c \in \mathbf{C}$ such that $\mathcal{K}(c) \subred w^*$. %(otherwise $\chi_{w^*} \in L_{\mathbf{C}}$, hence $w^*$ validates $\chi_{w^*}$, which is a contradiction: recall the properties of the splitting formula $\chi_V$, for a finite rooted transitive Kripke frame $V$, as stated in Proposition~\ref{prop:subred}).
Let $f \colon G \longrightarrow w^*$ be a surjective p-morphism, with $G \subseteq \mathcal{K}(c)$ generated subframe, witnessing $\mathcal{K}(c) \subred w^*$, and let $h \colon d \longrightarrow c$ in $G$ be such that $f(h) = w$. The restriction $h^* \longrightarrow w^*$ of $f$ is a surjective p-morphism, too (because $f(h^*) = f(h)^* = w^*$). Moreover, $h^*$ is the image of $\mathcal{K}(d)$ via $\mathcal{K}(h) \colon \mathcal{K}(d) \longrightarrow \mathcal{K}(c)$ and $\mathcal{K}(h)(\id_d) = h$. Composing $\mathcal{K}(d) \longrightarrow h^* \longrightarrow w^*$ with the inclusion of $w^*$ in $W$ yields a p-morphism $g \colon \mathcal{K}(d) \longrightarrow W$ such that $\varepsilon_W(g) = g(\id_d) = w$.
\endproof

\theorem\label{teo:presh}
Let $\mathbf{C}$ be the empty category or a one-object-category whose monoid of arrows is one of the following:
%The following list provides examples of categories $\mathbf{C}$ for which the hypotheses of \emph{Theorem~\ref{teo:suff}} are satisfied: $\mathbf{C}$ is the empty category ($L_{\mathbf{C}}$ is the inconsistent logic); $\mathbf{C}$ is a one-object-category whose monoid of arrows is one between
\begin{enumerate}
    \item $(\Z/2\Z,\times)$;% ($L_{\mathbf{C}} = S4.2 + \textit{bd}_2 + \textit{be}_1 + \textit{bi}_1$),
    \item $(\Z/3\Z,\times)$;% ($L_{\mathbf{C}} = S4.2 + \textit{bd}_2 + \textit{be}_1 + \textit{bi}_2$),
    \item the trivial monoid;% ($L_{\mathbf{C}} = S5 + \textit{be}_1$),
    \item $(\Z/2\Z,+)$.% ($L_{\mathbf{C}} = S5 + \textit{be}_2$).
\end{enumerate}
Then $\widehat{\mathbf{C}} \simeq L_\mathbf{C}\KFr_\lf$.
\endtheorem
\proof
As an example, consider the case 1. and call $\star$ the unique object of $\mathbf{C}$. We have that $\mathcal{K}(\star) \cong [2]$ and, for each Kripke frame $W$, the p-morphism $\varepsilon_W$ can be rewritten as follows
\[\begin{tikzcd}[row sep=tiny]
	{\KFr[[2],W]} & W \\
	{g \colon [2] \longrightarrow W} & {g(1).}
	\arrow["{\varepsilon_W}", from=1-1, to=1-2]
	\arrow[maps to, from=2-1, to=2-2]
\end{tikzcd}\]
We prove that $\varepsilon_W$ is injective for all $W \in \KFr$: given two p-morphisms $f, g \colon [2] \longrightarrow W$, if $f(1) = g(1)$, then the image
$$f([2]) = f(\up{1}) = \up{f(1)} = \up{g(1)} = g(\up{1}) = g([2])$$
is isomorphic to either $[1]$ or $[2]$, forcing $f(2) = g(2)$. % Moreover, the logic of $\mathcal{K}(\star) = [2]$ is exactly $S4.2 + \textit{bd}_2 + \textit{be}_1 + \textit{bi}_1$.
The claim then follows from Theorem~\ref{teo:suff}.
\endproof

We are now ready to prove that the $4$ classes listed in Theorem~\ref{teo:exact} (in addition to the class of locally finite frames for the inconsistent logic, namely $\{\emptyset\}$) are all presheaf toposes, hence Barr exact categories.
\theorem
Let $L$ be a normal modal logic with the finite model property and assume that $L \in \NExt \textup{\textsf{S4}}$. The following are equivalent:
\begin{enumerate}
    \item $L\KFr_\lf$ is Barr exact;
    \item $L\KFr_\lf$ is a Grothendieck topos;
    \item $L$ is the inconsistent logic, or $L\KFr_\lf$ is one of the following:
    \begin{enumerate}
        \item[(i)] $\{W \in \textup{\textsf{S4}}\KFr_\lf\ \vert\ \depth(W) \leq 2,\ \width(W) \leq 1,\ \exter(W) \leq 1,\ \inter(W) \leq 1\}$;
        \item[(ii)] $\{W \in \textup{\textsf{S4}}\KFr_\lf\ \vert\ \depth(W) \leq 2,\ \width(W) \leq 1,\ \exter(W) \leq 1,\ \inter(W) \leq 2\}$;
        \item[(iii)] $\{W \in \textup{\textsf{S4}}\KFr_\lf\ \vert\ \depth(W) \leq 1,\ \exter(W) \leq 1\}$;
        \item[(iv)] $\{W \in \textup{\textsf{S4}}\KFr_\lf\ \vert\ \depth(W) \leq 1,\ \exter(W) \leq 2\}$.
    \end{enumerate}
\end{enumerate}
\endtheorem
\proof
It remains to prove that the categories in (i)-(iv) are all Grothendieck toposes. As an example, the fact that the category in (i) is a Grothendieck topos follows from 1. of Theorem~\ref{teo:presh}: in that case, $\{\mathcal{K}(c)\ \vert\ c \in \mathbf{C}\} = \{[2]\}$ and it is straightforward to verify that the class in (i) coincides with $\{W \in \KFr_\lf\ \vert\ (\forall w \in W)([2] \subred w^*)\}$.
\endproof

Whenever the hypotheses of Theorem~\ref{teo:suff} are satisfied, $\mathcal{K}$ yields an equivalence between the category of presheaves $\widehat{\mathbf{C}}$ and $L_{\mathbf{C}}\KFr_\lf$, where $L_{\mathbf{C}}$ is the logic of $\{\mathcal{K}(c)\ \vert\ c \in \mathbf{C}\}$. In particular, $L_\mathbf{C}$ extends \textsf{S4}. %(as we observed, $\mathcal{K}(F)$ is a preorder for each $F \in \widehat{\mathbf{C}}$).
We can also say something for the irreflexive and transitive case, i.e.\ for those logics extending \textsf{GL}. Let us modify the construction of $\mathcal{K}$ in the following way: as before, consider the slice category $\mathbf{C}/F$, given a small category $\mathbf{C}$ and a presheaf $F \in \widehat{\mathbf{C}}$; this time, set
$$x \prec y \iff y = x h \text{ for some } h \colon d \longrightarrow c \text{ in } \mathbf{C} \text{ such that } \mathbf{C}[c,d] = \emptyset,$$
for $x \colon c \longrightarrow F$ and $y \colon d \longrightarrow F$ in $\mathbf{C}/F$. Let us call $\mathcal{K}'(F)$ the Kripke frame obtained in this way. Observe that $\mathcal{K}'(F)$ is irreflexive and transitive, hence $\mathcal{K}'(F) \in \textsf{GL}\KFr$. Define $\mathcal{K}'$ on morphisms as we did before.
\lemma
$\mathcal{K}'$ defines a faithful functor $\widehat{\mathbf{C}} \longrightarrow \KFr$.
\endlemma

As before, we define a functor $\mathcal{F}' \colon \KFr \longrightarrow \widehat{\mathbf{C}}$, by setting $\mathcal{F}'(W)(c) \coloneqq \KFr[\mathcal{K}'(c),W]$ for each $W \in \KFr$ and $c \in \mathbf{C}$. All the results that we proved for $\mathcal{K}$ and $\mathcal{F}$ can be adapted to $\mathcal{K}'$ and $\mathcal{F}'$ with only minor modifications (in particular, the proof of Theorem~\ref{teo:suff} requires the observation that the injectivity of $\varepsilon'_W$, counit of the adjunction $\mathcal{K}' \dashv \mathcal{F}'$, forces $\mathcal{K}'(c)$ to coincide with $\id_c^*$ for each $c \in \mathbf{C}$). %, except for Theorem~\ref{teo:suff}, which seems not to hold in some cases. However, if we assume furthermore that the category $\mathbf{C}$ itself is a poset, then Theorem~\ref{teo:suff} can be adapted as well.

Let $\mathbf{C}_\omega$ be the category whose objects are the positive natural numbers, having exactly one arrow $n \longrightarrow m$ whenever $n \leq m$. Moreover, let $\mathbf{C}_n$ be the full subcategory having $\{1,\dots,n\}$ as objects. For a positive natural number $m$, denote by $[m]'$ the $m$-elements linearly ordered irreflexive chain ($[m]' = \{1,\dots,m\}$, with $i \prec j$ if and only if $j < i$).
\theorem\label{teo:GL}
If $n$ is a positive natural number or $n = \omega$, then
$$\widehat{\mathbf{C}_n} \simeq \{W \in \textup{\textsf{GL.3}}\KFr_\lf\ \vert\ \depth(W) \leq n\}.$$
%$\mathcal{K}'$ restricts to an equivalence of categories between $\widehat{\mathbf{C}}$ and $L_{\mathbf{C}}\KFr_\lf$, where $L_{\mathbf{C}} = \textit{GL}.3 + \textit{bd}_n$ ($L_{\mathbf{C}}$ is the inconsistent logic if $n = 0 = \emptyset$).
\endtheorem
\proof
For each $m \in \mathbf{C}_n$, we have that $\mathcal{K}'(m) \cong [m]'$ and, for each Kripke frame $W$, the p-morphism $\varepsilon'_W$ can be rewritten as follows
\[\begin{tikzcd}[row sep=tiny]
	{\coprod_{m \in \mathbf{C}_n} \KFr[[m]',W]} & {W} \\
	{g \colon [m]' \longrightarrow W} & {g(m).}
	\arrow["{\varepsilon'_W}", from=1-1, to=1-2]
	\arrow[maps to, from=2-1, to=2-2]
\end{tikzcd}\]
It is straightforward to show that $\varepsilon'_W$ is injective for all $W \in \KFr$ (p-morphisms having an irreflexive chain $[m]'$ as domain are isomorphisms and the unique p-morphism $[m]' \longrightarrow [m]'$ is the identity). The claim now follows from (modified) Theorem~\ref{teo:suff}, observing that the class of locally finite frames generated by $\{[m]'\ \vert\ m \in \mathbf{C}_n\}$ is precisely the class of locally finite \textsf{GL.3}-frames of depth at most $n$.
\endproof

\section{Conclusions}
In this paper, motivated by the results of \cite{GdB24,infinitary}, we studied regularity and Barr exactness of $L\KFr_\lf$, the dual of the category $\Pro L\MA_\fin$ of profinite $L$-algebras, among the normal modal logics $L \supseteq \textsf{S4}$ with the finite model property. In doing so, we described the finite limits of $L\KFr_\lf$. Moreover, we introduced an adjunction between $\KFr$ and the category of presheaves over any small category $\mathbf{C}$ that allowed us to present some $L\KFr_\lf$ as presheaf toposes.

As we observed, the classification of the regularity of $L\KFr_\lf$ above \textsf{S4}, which is equivalent to an infinitary version of the Craig's Interpolation theorem, %which is equivalent to $L\KFr_\lf$ being regular and $\Pro L\MA_\fin$ and $L\MA_\fin$ satisfying the (AP),
subsumes the necessary conditions of Maksimova's parallel results for the finitary fragment \cite{Maksimova1979InterpolationTI,Maksimova1980InterpolationTI}. However, the complete characterizations must differ: as observed in Remark~\ref{rem:exc}, $\textsf{Grz.3}\MA_\fin$ satisfies the (AP), while $\textsf{Grz.3}\MA$ does not. In addition, regularity can be encountered in many other cases above \textsf{K4}, as we showed in Remark~\ref{rem:Horn}.

Passing to Barr exactness, the picture changes drastically: there are exactly five logics $L$ extending \textsf{S4} and such that $L\KFr_\lf$ is Barr exact. These five logics are all of finite height less than $3$
%rev locally finite
and of little interest (for such logics, $L\KFr_\lf$
%rev (the category of their locally finite Kripke frames
is equivalent to the topos of the actions of a finite small monoid). Still, one may ask whether Barr exactness is so rare everywhere: in particular a classification over \textsf{K4} cannot include only trivial examples, as witnessed by Theorem~\ref{teo:GL}.

Given that Barr exactness is rare, one may investigate almost Barr exactness, where a category is said to be \emph{almost Barr exact} if and only if it is regular and effective descent morphisms coincide with regular epimorphisms. Almost Barr exactness, in fact, turns out to be equivalent to effectiveness of special kinds of equivalence relations \cite{Beyond}, so that there is a concrete chance that it should hold  for some relevant categories of locally finite Kripke frames.

%\begin{references}
%\item A. Grothendieck (1957), Sur quelques points d'alg\'ebre homologique, T\^ohoku Math. J. \em{3}, 120-221.
%\end{references}

\bibliographystyle{plain}
\bibliography{references}

@book{misha,
    AUTHOR = {Chagrov, Alexander and Zakharyaschev, Michael},
     TITLE = {Modal logic},
    SERIES = {Oxford Logic Guides},
    VOLUME = {35},
      NOTE = {Oxford Science Publications},
 PUBLISHER = {The Clarendon Press, Oxford University Press, New York},
      YEAR = {1997},
     PAGES = {xvi+605},
      ISBN = {0-19-853779-4},
   MRCLASS = {03B45 (03-02)},
  MRNUMBER = {1464942},
MRREVIEWER = {J. M. Plotkin},
}

@article {GdB24,
    AUTHOR = {De Berardinis, Matteo and Ghilardi, Silvio},
     TITLE = {Profiniteness, monadicity and universal models in modal logic},
   JOURNAL = {Ann. Pure Appl. Logic},
  FJOURNAL = {Annals of Pure and Applied Logic},
    VOLUME = {175},
      YEAR = {2024},
    NUMBER = {7},
     PAGES = {Paper No. 103454, 25},
      ISSN = {0168-0072,1873-2461},
   MRCLASS = {03B45 (03G05 03G30 18C20)},
  MRNUMBER = {4736189},
MRREVIEWER = {R.\ Sh.\ Grigolia},
       DOI = {10.1016/j.apal.2024.103454},
       URL = {https://doi.org/10.1016/j.apal.2024.103454},
}

@book{stone,
  title={Stone spaces},
  author={Johnstone, Peter T},
  volume={3},
  year={1982},
  publisher={Cambridge University Press}
}

@book{GZ,
 author = {Ghilardi, Silvio and Zawadowski, Marek },
 title = {Sheaves, Games, and Model Completions: A Categorical Approach to Nonclassical Propositional Logics},
 year = {2011},
 isbn = {9048160367, 9789048160365},
 publisher = {Springer Publishing Company, Incorporated},
}

@article {leo,
    AUTHOR = {Esakia, Leo and Grigolia, Revaz},
     TITLE = {Christmas trees. {O}n free cyclic algebras in some varieties
              of closure algebras},
   JOURNAL = {Polish Acad. Sci. Inst. Philos. Sociol. Bull. Sect. Logic},
  FJOURNAL = {Polish Academy of Sciences. Institute of Philosophy and
              Sociology. Bulletin of the Section of Logic},
    VOLUME = {4},
      YEAR = {1975},
    NUMBER = {3},
     PAGES = {95--102},
      ISSN = {0138-0680},
}

@article {valentin,
  author = {Shehtman, Valentin B.},
  title = {Rieger-{N}ishimura lattices},
  journal = {Dokl. Akad. Nauk SSSR},
  volume = {19},
  pages = {1014--1018},
  year = {1978},
  }

@article {be1,
    AUTHOR = {Bellissima, Fabio},
     TITLE = {Finitely generated free {H}eyting algebras},
   JOURNAL = {J. Symbolic Logic},
  FJOURNAL = {The Journal of Symbolic Logic},
    VOLUME = {51},
      YEAR = {1986},
    NUMBER = {1},
     PAGES = {152--165},
      ISSN = {0022-4812},
   MRCLASS = {03G25 (03B55 03F55 06D20 08B20)},
  MRNUMBER = {830082},
MRREVIEWER = {R. Sh. Grigolia},
       DOI = {10.2307/2273952},
       URL = {https://doi.org/10.2307/2273952},
}

@article {be2,
    AUTHOR = {Bellissima, Fabio},
     TITLE = {An effective representation for finitely generated free
              interior algebras},
   JOURNAL = {Algebra Universalis},
  FJOURNAL = {Algebra Universalis},
    VOLUME = {20},
      YEAR = {1985},
    NUMBER = {3},
     PAGES = {302--317},
      ISSN = {0002-5240},
   MRCLASS = {08B20 (03G25 06D20)},
  MRNUMBER = {811691},
MRREVIEWER = {A. Goetz},
       DOI = {10.1007/BF01195140},
       URL = {https://doi.org/10.1007/BF01195140},
}

@incollection{ungh,
  title={Irreducible models and definable embeddings},
  author={Silvio Ghilardi},
  booktitle = {Logic Colloquium '92},
    series = {Stud. Logic Lang. Inform.},
     pages = {95--113},
 publisher = {CSLI Publ., Stanford, CA},
      year = {1995},
}

@article{Maksimova1979InterpolationTI,
  title={Interpolation theorems in modal logics and amalgamable varieties of topological Boolean algebras},
  author={Larisa L. Maksimova},
  journal={Algebra and Logic},
  year={1979},
  volume={18},
  pages={348-370},
  url={https://api.semanticscholar.org/CorpusID:121706349}
}

@article{Maksimova1982AbsenceOT,
  title={Absence of the interpolation property in the consistent normal modal extensions of the dummett logic},
  author={Larisa L. Maksimova},
  journal={Algebra and Logic},
  year={1982},
  volume={21},
  pages={460-463},
  url={https://api.semanticscholar.org/CorpusID:122109897}
}

@article{Maksimova1980InterpolationTI,
  title={Interpolation theorems in modal logics. {S}ufficient conditions},
  author={Larisa L. Maksimova},
  journal={Algebra and Logic},
  year={1980},
  volume={19},
  pages={120-132},
  url={https://api.semanticscholar.org/CorpusID:123207182}
}

@article{Maksimova1977CraigsTI,
  title={Craig's theorem in superintuitionistic logics and amalgamable varieties of pseudo-boolean algebras},
  author={Larisa L. Maksimova},
  journal={Algebra and Logic},
  year={1977},
  volume={16},
  pages={427-455},
  url={https://api.semanticscholar.org/CorpusID:120780596}
}

@article{Ghilardi1992QuantifiedEO,
  title={Quantified extensions of canonical propositional intermediate logics},
  author={Silvio Ghilardi},
  journal={Studia Logica},
  year={1992},
  volume={51},
  pages={195-214},
  url={https://api.semanticscholar.org/CorpusID:43177153}
}

@article {infinitary,
    AUTHOR = {De Berardinis, Matteo and Ghilardi, Silvio},
     TITLE = {A Proof Theory for Profinite Modal Algebras},
   JOURNAL = {ArXiv},
      YEAR = {2025},
    NOTE = {Submitted to a Springer volume in honor of M. Gehrke. \url{https://arxiv.org/abs/2507.06007/}},
}

@book {Tak,
    AUTHOR = {Takeuti, Gaisi},
     TITLE = {Proof theory},
    SERIES = {Studies in Logic and the Foundations of Mathematics},
    VOLUME = {Vol. 81},
 PUBLISHER = {North-Holland Publishing Co., Amsterdam-Oxford; American
              Elsevier Publishing Co., Inc., New York},
      YEAR = {1975},
     PAGES = {vi+372},
   MRCLASS = {02DXX (02-02)},
  MRNUMBER = {536648},
MRREVIEWER = {},
}

@article {thomason,
    AUTHOR = {Thomason, S. K.},
     TITLE = {Categories of frames for modal logic},
   JOURNAL = {J. Symbolic Logic},
  FJOURNAL = {The Journal of Symbolic Logic},
    VOLUME = {40},
      YEAR = {1975},
    NUMBER = {3},
     PAGES = {439--442},
      ISSN = {0022-4812},
   MRCLASS = {02C10 (02J05)},
  MRNUMBER = {381940},
MRREVIEWER = {A. Tauts},
       DOI = {10.2307/2272167},
       URL = {https://doi.org/10.2307/2272167},
}

@book{CWM,
  title={Categories for the working mathematician},
  author={Mac Lane, Saunders},
  volume={5},
  year={2013},
  publisher={Springer Science \& Business Media}
}

@incollection{manes,
  title={A triple theorethic construction of compact algebras},
  author={Ernest Manes},
  booktitle = {Seminar on Triples and Categorical Homology Theory},
     pages = {91--118},
 publisher = {Springer},
      year = {1969},
}

@article{gehrke,
title = {Stone duality, topological algebra, and recognition},
journal = {Journal of Pure and Applied Algebra},
volume = {220},
number = {7},
pages = {2711-2747},
year = {2016},
issn = {0022-4049},
doi = {https://doi.org/10.1016/j.jpaa.2015.12.007},
author = {Mai Gehrke},
}

@book {MR2522659,
    AUTHOR = {Lurie, Jacob},
     TITLE = {Higher topos theory},
    SERIES = {Annals of Mathematics Studies},
    VOLUME = {170},
 PUBLISHER = {Princeton University Press, Princeton, NJ},
      YEAR = {2009},
     PAGES = {xviii+925},
      ISBN = {978-0-691-14049-0; 0-691-14049-9},
   MRCLASS = {18-02 (18B25 18E35 18G30 18G55 55U40)},
  MRNUMBER = {2522659},
MRREVIEWER = {Mark\ Hovey},
       DOI = {10.1515/9781400830558},
       URL = {https://doi.org/10.1515/9781400830558},
}

@incollection {Beyond,
    AUTHOR = {Janelidze, George and Sobral, Manuela and Tholen, Walter},
     TITLE = {Beyond {B}arr exactness: effective descent morphisms},
 BOOKTITLE = {Categorical foundations},
    SERIES = {Encyclopedia Math. Appl.},
    VOLUME = {97},
     PAGES = {359--405},
 PUBLISHER = {Cambridge Univ. Press, Cambridge},
      YEAR = {2004},
      ISBN = {0-521-83414-7},
   MRCLASS = {18E10},
  MRNUMBER = {2056587},
}

@book {MR1953060,
    AUTHOR = {Johnstone, Peter T.},
     TITLE = {Sketches of an elephant: a topos theory compendium. {V}ol. 1},
    SERIES = {Oxford Logic Guides},
    VOLUME = {43},
 PUBLISHER = {The Clarendon Press, Oxford University Press, New York},
      YEAR = {2002},
     PAGES = {xxii+468+71},
      ISBN = {0-19-853425-6},
   MRCLASS = {18B25 (18-02)},
  MRNUMBER = {1953060},
MRREVIEWER = {Colin\ McLarty},
}

@book {MR163850,
    AUTHOR = {Rasiowa, Helena and Sikorski, Roman},
     TITLE = {The mathematics of metamathematics},
    SERIES = {Monografie Matematyczne [Mathematical Monographs]},
    VOLUME = {Tom 41},
 PUBLISHER = {Pa\'nstwowe Wydawnictwo Naukowe, Warsaw},
      YEAR = {1963},
     PAGES = {522},
   MRCLASS = {02.52},
  MRNUMBER = {163850},
MRREVIEWER = {G.\ Kreisel},
}

@article {MR1171642,
    AUTHOR = {Chagrov, Alexander and Zakharyashchev, Michael},
     TITLE = {Modal companions of intermediate propositional logics},
   JOURNAL = {Studia Logica},
  FJOURNAL = {Studia Logica. An International Journal for Symbolic Logic},
    VOLUME = {51},
      YEAR = {1992},
    NUMBER = {1},
     PAGES = {49--82},
      ISSN = {0039-3215,1572-8730},
   MRCLASS = {03B45 (03B20 03B55)},
  MRNUMBER = {1171642},
MRREVIEWER = {Stanis\l aw\ J.\ Surma},
       DOI = {10.1007/BF00370331},
       URL = {https://doi.org/10.1007/BF00370331},
}

@book {MR3967703,
    AUTHOR = {Esakia, Leo},
     TITLE = {Heyting algebras},
    SERIES = {Trends in Logic---Studia Logica Library},
    VOLUME = {50},
    EDITOR = {Bezhanishvili, Guram and Holliday, Wesley H.},
   EDITION = {Russian},
      NOTE = {Duality theory},
 PUBLISHER = {Springer, Cham},
      YEAR = {2019},
     PAGES = {xv+95},
      ISBN = {978-3-030-12095-5; 978-3-030-12098-6; 978-3-030-12096-2},
   MRCLASS = {06-02 (06D20 06D50)},
  MRNUMBER = {3967703},
MRREVIEWER = {Marcel\ Ern\'e},
       DOI = {10.1007/978-3-030-12096-2},
       URL = {https://doi.org/10.1007/978-3-030-12096-2},
}

@article{Godel,
author = {Carai, Luca},
year = {2026},
month = {02},
pages = {1-24},
title = {FREE ALGEBRAS AND COPRODUCTS IN VARIETIES OF {G}{\"O}DEL ALGEBRAS},
journal = {The Journal of Symbolic Logic},
doi = {10.1017/jsl.2026.10194}
}

@book {MR1837791,
    AUTHOR = {Blackburn, Patrick and de Rijke, Maarten and Venema, Yde},
     TITLE = {Modal logic},
    SERIES = {Cambridge Tracts in Theoretical Computer Science},
    VOLUME = {53},
 PUBLISHER = {Cambridge University Press, Cambridge},
      YEAR = {2001},
     PAGES = {xxii+554},
      ISBN = {0-521-80200-8; 0-521-52714-7},
   MRCLASS = {03-02 (03B45 08A70 68T27)},
  MRNUMBER = {1837791},
MRREVIEWER = {J.\ M.\ Plotkin},
       DOI = {10.1017/CBO9781107050884},
       URL = {https://doi.org/10.1017/CBO9781107050884},
}
\end{document}